\documentclass[10pt,a4paper]{article}
\usepackage{preamble}
\usepackage[toc, page]{appendix}
\newtheorem*{thm*}{Theorem 4.10}
\newtheorem*{defn*}{Definition B.1}

\newcommand{\exedgename}{exceptional edge}
\newcommand{\Exedgename}{Exceptional edge}
\newcommand{\exedge}{\eta}
\newcommand{\exloopname}{exceptional loop}
\newcommand{\Exloopname}{Exceptional loop}
\newcommand{\exloop}{\ell}

\newcommand{\hue}{hue}

\newcommand{\hues}{hues}

\newcommand{\hueOrbit}{colour}

\newcommand{\hueOrbits}{colours}

\author{Michelle Strumila \\ \textit{Monash University}\footnote{Most of this work was completed during the author's PhD at the University of Melbourne, under Dr Marcy Robertson}}
\title{Quasi Modular Operads}
\begin{document}
\maketitle
\begin{abstract}
    Modular operads are an extension of operads.  In the same way that operads, as dendroidal sets, can be considered as presheaves over the category of trees, so can modular operads be considered as presheaves over a category of graphs, leading towards infinity modular operads.
    
    This paper contains a definition of the Kan condition for infinity modular operads (Definition \ref{def:kancond}), as well as a proof of the Nerve Theorem (Theorem \ref{thm:GraphicalNerve}) and the equivalence of the Kan and Segal conditions.  Appendix \ref{app:cycOp} contains the same material for cyclic operads.
\end{abstract}
\tableofcontents
\section{Introduction}
Modular operads were introduced by \citet{getzler1998modular} to aid in Feynman diagram calculations, and were further studied by \citet{robertson2019modular}, and \citet{raynor2019distributivelaw}.  They can be thought of as operads of higher genus, and this higher genus nature makes them useful for studying surfaces \citep{fong2015decorated, giansiracusa2011framed, MR1752293, Segal1988, teleman2016five, tillmann2000higher}.

Infinity categories and infinity operads have been previously defined under various models \citep{bousfield1972homotopy, eilenberg1950semi, moerdijk2007, moerdijk2009inner}, with the aid of simplicial and dendroidal sets.  This paper, along with \citep{Robertson2019ModNerve} and \citep{robertson2019modular}, provides a higher genus analogue to simplicial and dendroidal sets.

The following table illustrates the relationships between categories, modular operads, and other generalised operads, as well as their underlying graphs.  Just as morphism composition in categories is linear, and operads are tree-like, modular operads are graph-like.  An analysis of these various categories of graphs is given by \citet{Hackney2024}.
	\begin{center}
	\begin{tabular}{|c|ccc|} 
	\hline 
	 & Paths & Trees & Graphs \\ 
	\hline 
	Directed & Categories & Operads & Wheeled properads  \\ 
	Undirected & Dagger categories & Cyclic operads & Modular operads \\ 
	\hline 
	\end{tabular} 
	\end{center}
Wheeled properads are similar to modular operads, but use the category of directed graphs.  They are properads in which the morphisms may have loops, with properads being operads in which the morphisms may have both multiple inputs and multiple outputs.  More information on properads and wheeled properads can be found in \citep{hackney2014infinity}.  

There are different models for infinity categories.  \citet{bergner2010survey} gives equivalences between quasi categories, enriched categories, Segal categories, 
and complete Segal spaces.  Of particular note are quasi categories and Segal categories, as these are both defined as simplicial sets satisfying a particular condition, respectively the weak inner Kan condition and the Segal condition.  Similarly, there are different ways of defining infinity modular operads.  

The Nerve Theorem for categories states the relationship between simplicial sets and quasi categories \citep{grothendieck1960technique}.  For modular operads, \citet{robertson2019modular, Robertson2019ModNerve} provide definitions of graphical sets and the modular Segal condition, and a relationship between them.

This paper extends some of their work on the category of graphs and provides the definition of the weak inner Kan condition for modular operads. 
Theorem \ref{thm:modularsegalkan} shows the equivalence between the Segal and Kan conditions, and thus completes the nerve theorem for modular operads, given fully as Theorem \ref{thm:GraphicalNerve}.

\begingroup
\def\thetheorem{\ref{thm:GraphicalNerve}}
\begin{thm*}[Graphical Nerve Theorem]
\sloppy
	Let $X$ be a graphical set ${X:\Graph^{op} \to \cat{Set}}$.  Then the following are equivalent:
	\begin{enumerate}
		\item There exists a modular operad $\M$ such that $n_g(\M) = X$
		\item $X$ satisfies the (strict) inner Kan condition
		\item $X$ satisfies the Segal condition
	\end{enumerate}
\end{thm*}
\addtocounter{thm}{-1}
\endgroup

\section{Modular Operads}
There are a number of definitions given for Modular operads in the literature, both coloured and monochromatic.  Monochromatic operads were first introduced by \citet{getzler1998modular}.

Coloured modular operads, defined below, have only recently been studied in the development of a theory of infinity modular operads.  Definition \ref{def:modop} is based on that of \citet{raynor2019distributivelaw}, and similar to \citet{robertson2019modular} and \citet[Definition 3.3.1]{hinich2002cyclic}.  An equivalent definition is given by \citet{Stoll2022}, where modular operads are defined as an algebra over a certain kind of properad.  
Elsewhere in the literature, coloured modular operads are sometimes called Compact Symmetric Multicategories
\citep{joyal2011feynman}.

Note that the adjective ``coloured'' will often be dropped. Throughout this paper, we let $\E = (\E, \otimes, 1)$ be a closed symmetric monoidal category; actually, we only require $\E$ to be cartesian.  The main categories used in this paper are $\cat{Set}$, the category of sets and functions, and $\cat{Grp}$, the category of groups and group homomorphisms.

	\begin{defn}[Profile]
	\label{def:profile}
	Given a set $A$, a profile $\underline{c}$ over $A$ of length $n$ is a list of elements of $A$, $$\underline{c} = (c_1, c_2, \ldots c_n)$$
\end{defn}
\begin{defn}[Sub-profile]
	Given a profile $\underline{c} = (c_1, c_2, \ldots c_n)$, a sub-profile $\underline{d}$ of $\underline{c}$ is a list of elements of $\underline{c}$,  listed in $\underline{d}$ in the same order that they appear in $\underline{c}$. 
\end{defn}
\begin{defn}[Modular operad] 
	\label{def:modop}
	A coloured modular operad, $\mathcal{M}$, over a closed symmetric monoidal category $\mathcal{E}$, is defined by the following.
	\begin{itemize}
		\item There is a set of objects $ob(\mathcal{M})$, sometimes called \hues, equipped with an involution $\iota$.  Each orbit is known as a \hueOrbit.
		\item For each profile $\underline{c}$ over $ob(\mathcal{M})$, there is an object of $\mathcal{E}$, $\mathcal{M} (\underline{c})$.  Usually $\mathcal{M} (\underline{c})$ will be a set.
		\item There is a right action of the symmetric group.  That is, for any profile $\underline{c} = (c_{i_1}, \ldots c_{i_n} )$
		 and permutation $\sigma \in \Sigma_n$, there is a bijection $\M(c_1, \ldots, c_n) \to \M (c_{\sigma(1)}, \ldots, c_{\sigma(n)})$
		\item There is a composition.  Let $\underline{X} = (c_{i_1}, \ldots c_{i_r} )$ be a sub-profile of $\underline{c}$ such that $\left(\iota(c_{i_1}), \ldots \iota(c_{i_r}) \right) = \left( d_{j_1}, \ldots d_{j_r}\right)$ is a sub-profile of $\underline{d}$.   
		Then for $\theta \in \mathcal{M}(\underline{c})$ and $\theta' \in \mathcal{M}(\underline{d})$,
		$$\theta \circ_X \theta' \mapsto \theta''$$
		where $\theta'' \in \mathcal{M}(c_1, \ldots \hat{c_{i_1}}, \ldots \hat{c_{i_r}}, \ldots c_m, d_1, \ldots \hat{d_{j_1}}, \ldots \hat{d_{j_r}}, \ldots d_n )$. 
		Note that $\hat{c_{i_j}}$ is used to denote removing $c_{i_j}$, and that the removed $c_{i_j}$'s may be interlaced with the non-removed objects.  
		For ease of reading, let $$\underline{cd}_X = c_1, \ldots \hat{c_{i_1}}, \ldots \hat{c_{i_r}}, \ldots c_m, d_1, \ldots \hat{d_{j_1}}, \ldots \hat{d_{j_r}}, \ldots d_n.$$ 
		\item There is a contraction.  Let $\underline{c}$ be a profile in which there are at least two objects $c_i$ and $c_j$ such that $\iota(c_i) = c_j$.  Then there is a contraction, $\zeta_{ij}: \mathcal{M}(\underline{c}) \to \mathcal{M}(c_1, \ldots, \hat{c_i}, \ldots, \hat{c_j}, \ldots, c_n)$.
	\end{itemize}
	Moreover, these satisfy the following axioms:
	\begin{enumerate}
		\item Composition is associative.  That is, for all profiles $\underline{a}$, $\underline{b}$, and $\underline{c}$:
			\begin{center}
			\begin{tikzcd}
				\mathcal{M}(\underline{a})\otimes \mathcal{M}(\underline{b})\otimes \mathcal{M}(\underline{c}) \arrow{r}{\circ_Y} \arrow{d}{\circ_X} & \mathcal{M}(\underline{a}) \otimes \mathcal{M}(\underline{bc}_Y)\arrow{d}{\circ_{X}} \\
				\mathcal{M}(\underline{ab}_X) \otimes \mathcal{M}(\underline{c}) \arrow{r}{\circ_Y} & \mathcal{M}(\underline{abc}_{XY})
			\end{tikzcd}
			\end{center}
		\item Composition is unital.  That is, for each \hue{} $h$ there exists an identity $\eta_h$, such that, for all profiles $\underline{c}$ containing $\iota(h)$ and for all $\theta \in \mathcal{M}(\underline{c})$, $$\theta \circ \eta_h = \theta = \eta_h \circ \theta$$  In each case the composition occurs over the subprofile of $\theta$ consisting of the single element $\iota(h)$.
		\item Composition is equivariant.  That is, it commutes with the action of the symmetric group.  So for any two morphisms $\alpha \in \mathcal{M}(\underline{c})$ and $\beta \in \mathcal{M}(\underline{d})$, where $\underline{c}$ is of length $m$ and $\underline{d}$ is of length $n$, and any $\sigma \in \Sigma_n, \tau \in \Sigma_m$, 
		$$\alpha \circ_{\sigma(X)} \sigma(\beta) = \sigma'(\alpha \circ_X \beta )$$
		$$\tau(\alpha) \circ_X \beta = \tau'( \alpha \circ_i \beta)$$
		Where $\sigma' \in \Sigma_{n+m-2|u|}$ refers to the element that acts on $\underline{cd}_X$ by doing $\sigma$ on each $d_k$ and the identity on each $c_j$, and $\tau' \in \Sigma_{n+m-2|u|}$ does the identity on each $d_k$ and permutes each $c_j$ according to $\tau$.
		\item If $\underline{c}$ is a profile in which at least two contractions are possible, so $\iota(c_i) = c_j$ and $\iota(c_k) = c_\ell$, then:
			\begin{center}
			\begin{tikzcd}
				\mathcal{M}(\underline{c}) \arrow{r}{\zeta_{ij}} \arrow{d}{\zeta_{k\ell}} & \mathcal{M}(c_1, \ldots, \hat{c_i}, \ldots, \hat{c_j}, \ldots, c_n) \arrow{d}{\zeta_{k\ell}} \\
				\mathcal{M}(c_1, \ldots, \hat{c_k}, \ldots, \hat{c_{\ell}}, \ldots, c_n) \arrow{r}{\zeta_{ij}} & \mathcal{M}(c_1, \ldots, \hat{c_i}, \ldots, \hat{c_j}, \ldots, \hat{c_k}, \ldots, \hat{c_{\ell}},\ldots, c_n)
			\end{tikzcd}
			\end{center}
			Note that the profile $(c_1, \ldots, \hat{c_i}, \ldots, \hat{c_j}, \ldots, c_n)$ may equivalently be referred to as $\underline{c} \setminus c_i, c_j$.
		\item The contraction commutes with the composition.  That is, if $\underline{c}$ is a profile such that some $\iota(c_j) = c_k$, and neither $c_j$ nor $c_k$ is in $X$, then the following diagram commutes.
			\begin{center}
			\begin{tikzcd}
				\mathcal{M}(\underline{c}) \otimes \mathcal{M}(\underline{d}) \arrow{r}{\zeta_{jk}} \arrow{d}{\circ_X} & \mathcal{M}(\underline{c}\setminus c_j, c_k) \otimes \mathcal{M}(\underline{d}) \arrow{d}{\circ_X} \\
				\mathcal{M}(\underline{cd}_X) \arrow{r}{\zeta_{jk}} & \mathcal{M}(\underline{cd}_X \setminus c_j, c_k)
			\end{tikzcd}
			\end{center}
		\item Parallel gluing of distinct elements. Let $\underline{c}$ and $\underline{d}$ be two profiles, and let $\underline{X} = (c_{i_1}, \ldots c_{i_r} )$ be a sub-profile of $\underline{c}$ such that $\left(\iota(c_{i_1}), \ldots \iota(c_{i_r}) \right) = \left( d_{j_1}, \ldots d_{j_r}\right)$ is a sub-profile of $\underline{d}$.  Let $a$ and $b$ be two colours in $X$.  Then the following diagram commutes.
			\begin{center}
			\begin{tikzcd}
				\mathcal{M}(\underline{c}) \otimes \mathcal{M}(\underline{d}) \arrow{r}{\circ_{X \setminus a}} \arrow{d}{\circ_{X \setminus b}} & \mathcal{M}(\underline{cd}_{X \setminus a}) \arrow{d}{\zeta_{a}} \\
				\mathcal{M}(\underline{cd}_{X \setminus b}) \arrow{r}{\zeta_{b}} & \mathcal{M}(\underline{cd}_{X}) 
			\end{tikzcd}
			\end{center}
	\end{enumerate}
\end{defn}
Just as operations in an operad can be pictured as corollae, with composition by glueing into trees, operations in a modular operad can be pictured as corollae with loops.  Composition involves glueing edges together, but unlike operads, these operations may be glued into higher genus graphs, rather than just trees.  To compose two operations, glue some legs of one corolla to legs of the other which share the same \hueOrbit{}.  The contraction operation represents gluing together two legs of the same \hueOrbit{} attached to the same operation.

\begin{figure}
\label{fig:defAxioms}
\begin{center}
\includegraphics[width=0.45\textwidth]{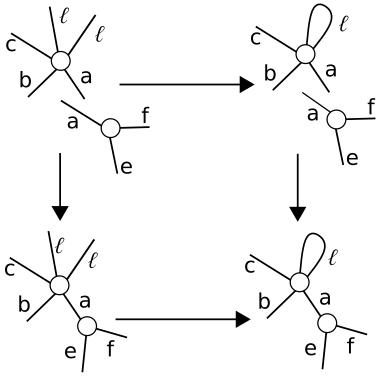} 
\includegraphics[width=0.45\textwidth]{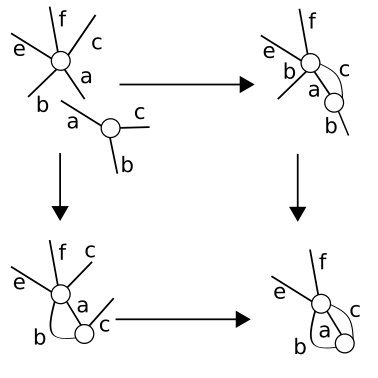} 
\caption{An illustration of axioms 5 (left) and 6 (right).}
\end{center}
\end{figure}

Note that Axioms 5 and 6 of Definition \ref{def:modop} are both referring to the composition commuting with contraction (see Figure \ref{fig:defAxioms}.  In Axiom 5, it is a local contraction that is not influenced by the composition, whereas in Axiom 6 it is two different contractions that become essentially equivalent by the composition which occurs.

Morphisms between modular operads can be defined in the same way as for other operadic generalisations, and then the category of modular operads can be defined.

\begin{defn}[Category of Modular operads, $\cat{ModOpd}$]
	A morphism $F: \mathcal{M} \to \mathcal{N}$ consists of the following information.
	\begin{itemize}
		\item A function $F: ob(\mathcal{M}) \to ob(\mathcal{N})$ which respects the involution.
		\item For each profile $\underline{c}$, a function $$F_{\underline{c}} : \Op(c_1, c_2, \ldots c_n) \to \oP_g(F(c_1), F(c_2), \ldots F(c_n))$$
			 such that identities, symmetric group actions, composition, and contraction are all preserved.
	\end{itemize}
	Then this is the category of modular operads $\cat{ModOpd}$
\end{defn}

\begin{eg}[Monochromatic modular operad]
	When the set $ob(\mathcal{M})$ contains only one object, the definition reduces to that of a monochromatic modular operad.
\end{eg}
\section{Graphical Sets and the Segal condition}
	\subsection{Graphical category}
			There are multiple ways of defining a graph.  Definition \ref{def:graph} is used so as to match the usual definition in the modular operads literature.  For an examination of the various definitions and proofs of their equivalence, see  \citet[Proposition 15.8]{batanin2017homotopy}.  
	The following definitions, leading up to and including $\Graph$, can be found in \citet[Section 1]{robertson2019modular}.  
    There is a description of face and degeneracy maps, as well as a table of relations, in Appendix \ref{sec:facemaprelations}.

\begin{defn}[Graph]
	\label{def:graph}
	A graph consists of sets $E, H, V \subseteq \mathcal{F}$ and  maps $i$, $s$, and $t$, satisfying the diagram 
	\begin{center}
	\begin{tikzcd}
		E \arrow[loop left]{l}{i} & H \arrow[swap]{l}{s} \arrow{r}{t} & V \\
	\end{tikzcd}
	\end{center}
	where $\mathcal{F}$ is some infinite set, $s$ is a monomorphism, and $i$ is a fixed point free involution.
\end{defn}
For those who think in terms of a more topological definition of graph, the letters these sets are denoted by above provide some clue as to the relationship between these definitions.  The set $V$ is the set of vertices, and the set $E$ is the set of arcs, or half edges, with an involution sending half edges to their opposite direction counterparts.  An edge would then be defined as a set $\{a,i(a)\}$, where $a$ is some half edge.

The set $H$ contains all those half edges which point towards a vertex.  The map $t$ sends each to its target vertex, which the map $s$ includes into the set of all arcs, $E$.  

\begin{defn}[Internal and external edges]
\label{def:internalEdge}
	The internal edges of a graph $G$ are those edges which have a vertex on each side (possibly the same vertex).  That is, they are those edges of the form $\{sd, sd'\}$, where $d,d'\in H$. Edges with at least one `loose end', i.e. of the form $\{sd, a\}$ where $d \in H$ and $a \notin H$, are called external edges, or legs.
\end{defn}
\begin{defn}[Boundary of a graph]
\label{def:boundary}
	The boundary $\partial (G)$ of a graph $G$ is the set $E \setminus sH$.  That is, the set of all edges where (at least) one side is unattached to a vertex.
\end{defn}
\begin{defn}[Neighbourhood]
\label{def:neighbourhood}
	Given a graph $G$ and a vertex $v$, the neighbourhood $nb(v)$ of $v$ is the set of all legs of $v$.  That is, the edges associated with the preimage $i \circ s \circ t^{-1}(v) \subset E$.
\end{defn}

Note that Definition \ref{def:graph} does not adequately capture the existence of the graph $S^1$, consisting of a single edge looped into itself, and no vertices.  Or rather, it cannot distinguish between a single edge with no vertices and the loop with no vertices.  Thus, we define them both here and append to the definition of graph.

\begin{defn}[\Exedgename]
\label{def:exceptionalEdge}
The \exedgename, $\exedge$, is the single isolated edge which is not connected to any vertices, topologically equivalent to the open interval $(0,1)$.  In the manner of Definition \ref{def:graph}, its sets are:
	\begin{center}
	\begin{tikzcd}
		\{a,a'\} \arrow[loop left]{l}{i} & \emptyset \arrow[swap]{l}{s} \arrow{r}{t} & \emptyset \\
	\end{tikzcd}
	\end{center}
	Here $i$ maps $a$ to $a'$ and vice versa.
\end{defn}

\begin{defn}[\Exloopname]
\label{def:exceptionalLoop}
The \exloopname, $\exloop$ is the single isolated edge, not connected to any vertices, which loops back on itself (topologically equivalent to a circle $S^1$).  In the manner of Definition \ref{def:graph}, its sets would be the same as the \exedgename:
	\begin{center}
	\begin{tikzcd}
		\{a,a'\} \arrow[loop left]{l}{i} & \emptyset \arrow[swap]{l}{s} \arrow{r}{t} & \emptyset \\
	\end{tikzcd}
	\end{center}
	where $i$ maps $a$ to $a'$ and vice versa.  However, unlike the \exedgename{} $\exedge$, the \exloopname{} $\exloop$ should not have any legs in the boundary.  Therefore we append $\exloop$ to the set of graphs, and define it to be an internal edge.
\end{defn}

From now on assume all graphs are connected.

\begin{defn}[Connected Graph]
\label{def:connected}
	Consider a graph $G$, 
	\begin{center}
	\begin{tikzcd}
		E \arrow[loop left]{l}{i} & H \arrow[swap]{l}{s} \arrow{r}{t} & V \\
	\end{tikzcd}
	\end{center}
	Then $G$ is called a \emph{connected graph} if at least one of the following is true:
	\begin{itemize}
		\item $G$ is the \exedgename{} $\exedge$ (Definition \ref{def:exceptionalEdge}).
		\item $G$ is the \exloopname{} $\exloop$ (Definition \ref{def:exceptionalLoop})
		\item $G$ satisfies all of the following:
		\begin{enumerate}
			\item $G$ is not the empty graph.
			\item All edges are attached to at least one vertex.  That is, the map $s: H \to E$ is surjective.  
			\item For any two distinct vertices $u$ and $v$, there exists at least one path between them.  That is, a list of edges $(e_j)_{j=1}^r$ and a list of vertices $(v_i)_{i=0}^{r}$ such that each vertex and edge is distinct, each $e_j$ is adjacent to both $v_{j-1}$ and $v_j$, and $u=v_0$ and $v=v_r$.
		\end{enumerate}
	\end{itemize}
\end{defn}

\begin{figure}
	\label{fig:barbell}
	\begin{center}
		\includegraphics[scale=1]{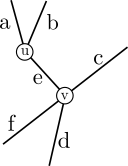}
	\end{center}
	\caption{A barbell graph, consisting of two vertices $u$ and $v$, an edge between them $e$, and some number of legs attached to each vertex.}
\end{figure}

One important graph is known as the \emph{barbell} graph (see Figure \ref{fig:barbell}).  This consists of two vertices, connected by an edge, with some number of edges emanating from them.  More formally:
\begin{defn}[Barbell]
    \label{def:barbell}
    A barbell is a graph consisting of the following sets and functions:
    \begin{align*}
        V & = \{u,v\} \\
        E & = \{e, e', u_1, u_2, \ldots, u_i, u'_1, u'_2, \ldots, u'_i, v_1, v_2, \ldots, v_j, v'_1, v'_2, \ldots, v'_j \} \\ 
        H & = \{e, e', u_1, u_2, \ldots, u_i, v_1, v_2, \ldots, v_j\} \\
        s(x) & = x \\
        t(x_i) & = x, \; t(e)=u, \; t(e')=v\\
        i(x) & = x', \; i(x') = x\\
    \end{align*}
    That is, the two vertices are labelled $u$ and $v$, with an edge $e$ between them, and edges $u_k$ and $v_k$ emanating from them.
\end{defn}

Graphical maps, which form the morphisms between graphs, shall be defined shortly.  
\begin{samepage}
In the following definitions let $G$ and $G'$ be graphs, where $G$ is 
	\begin{center}
	\begin{tikzcd}
		E \arrow[loop left]{l}{i} & H \arrow[swap]{l}{s} \arrow{r}{t} & V \\
	\end{tikzcd}
	\end{center}
	and $G'$ is
		\begin{center}
	\begin{tikzcd}
		E' \arrow[loop left]{l}{i'} & H' \arrow[swap]{l}{s'} \arrow{r}{t'} & V' \\
	\end{tikzcd}
	\end{center}
\end{samepage}

\begin{defn}[\'Etale]
	A natural transformation $f: G \to G'$ is an \'etale map if the right hand square is a pullback:
	\begin{center}
	\begin{tikzcd}
		E \arrow[loop left]{l}{i} \arrow{d}{f_E} & H \arrow[swap]{l}{s} \arrow{d}{f_H} \arrow{r}{t} & V \arrow{d}{f_V} \\
		E' \arrow[loop left]{l}{i'} & H' \arrow[swap]{l}{s'} \arrow{r}{t'} & V' \\
	\end{tikzcd}
	\end{center}
\end{defn}
Note that, when $G$ is the \exedgename{} or \exloopname{}, the map is \'etale because $H$ and $V$ are empty sets.  However, when $G'$ is the \exedgename{} or \exloopname{}, the map cannot be \'etale for the same reason.

\begin{eg}
	\label{eg:notetale}
	Consider two graphs, $G$ and $G'$, with a map $\varphi$ between them.   
	\begin{center}
		\includegraphics[scale=0.6]{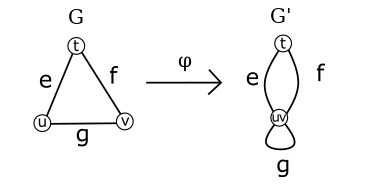}
	\end{center}
	Firstly $G$ is defined as:
	\begin{center}
	\begin{tikzcd}
		\{e, e^\dag, f, f^\dag, g, g^\dag \} \arrow[loop left]{l}{i} & \{e, e^\dag, f, f^\dag, g, g^\dag \} \arrow[swap]{l}{s} \arrow{r}{t} & \{t,u,v\} \\
	\end{tikzcd}
	\end{center}
	with the involution defined in the obvious way (i.e. $x \mapsto x^\dag$ and vice versa), $s$ defined by the inclusion, 
	and $t$ defined by 
	\begin{align*}
		e & \mapsto t & e^\dag & \mapsto u \\
		f & \mapsto v & f^\dag & \mapsto t \\
		g & \mapsto u & g^\dag & \mapsto v. \\
	\end{align*}
	Secondly $G'$ is defined as:
	\begin{center}
	\begin{tikzcd}
		\{e, e^\dag, f, f^\dag, g, g^\dag \} \arrow[loop left]{l}{i} & \{e, e^\dag, f, f^\dag, g, g^\dag \} \arrow[swap]{l}{s} \arrow{r}{t} & \{t,uv\} \\
	\end{tikzcd}
	\end{center}
	with the involution defined in the obvious way (i.e. $x \mapsto x^\dag$ and vice versa), $s'$ defined by the inclusion,
	and $t'$ defined by 
	\begin{align*}
		e & \mapsto t & e^\dag & \mapsto uv \\
		f & \mapsto uv & f^\dag & \mapsto t \\
		g & \mapsto uv & g^\dag & \mapsto uv. \\
	\end{align*}
	Then define the map $\varphi: G \to G'$, where $\varphi_E(x) \mapsto x$, $\varphi_H(x) \mapsto x$, and for the vertices $\varphi_V(t) \mapsto t$, $\varphi_V(u) \mapsto uv$, and $\varphi_V(v) \mapsto uv$.  Then this is not \'etale, because the commutative diagram
	\begin{center}
	\begin{tikzcd}
		H \arrow{d}{\varphi_H} \arrow{r}{t} & V \arrow{d}{\varphi_V} \\
		H' \arrow{r}{t'} & V' \\
	\end{tikzcd}
	\end{center}
	is not a pullback.
\end{eg}

\begin{defn}[Embedding]
	An embedding $G \to G'$ is an \'etale map for which $V \to V'$ is a monomorphism. 
\end{defn}

\begin{eg}[\'Etale but not embedding]
	\label{eg:etalenotemb}
	Consider two graphs, $G$ and $G'$, with a map $\varphi$ between them.   
	\begin{center}
		\includegraphics[scale=0.6]{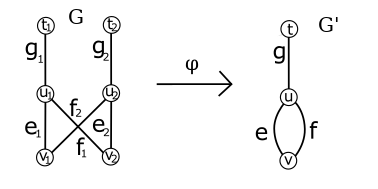}
	\end{center}
	Firstly $G$ is defined as:
	\begin{center}
	\begin{tikzcd}
		\{e_1, e_1^\dag, e_2, \ldots g_2^\dag \} \arrow[loop left]{l}{i} & \{e_1, e_1^\dag, e_2, \ldots g_2^\dag \} \arrow[swap]{l}{s} \arrow{r}{t} & \{t_1, t_2, u_1, u_2, v_1, v_2\} \\
	\end{tikzcd}
	\end{center}
	with the involution defined in the obvious way (i.e. $x \mapsto x^\dag$ and vice versa), $s$ the inclusion, 
	and $t$ defined by 
	\begin{align*}
		e_i & \mapsto u_i & e_i^\dag & \mapsto v_i \\
		f_i & \mapsto u_i & f_1^\dag & \mapsto v_2 \\
			& 			  & f_2^\dag & \mapsto v_1 \\
		g_i & \mapsto t_i & g_i^\dag & \mapsto u_i \\
	\end{align*}
	Secondly $G'$ is defined as:
	\begin{center}
	\begin{tikzcd}
		\{e, e^\dag, f, f^\dag, g, g^\dag \} \arrow[loop left]{l}{i} & \{e, e^\dag, f, f^\dag, g, g^\dag \} \arrow[swap]{l}{s} \arrow{r}{t} & \{t,u,v\} \\
	\end{tikzcd}
	\end{center}
	with the involution defined in the obvious way (i.e. $x \mapsto x^\dag$ and vice versa), $s'$ the inclusion,
	and $t'$ defined by 
	\begin{align*}
		e & \mapsto u & e^\dag & \mapsto v \\
		f & \mapsto u & f^\dag & \mapsto v \\
		g & \mapsto t & g^\dag & \mapsto u \\
	\end{align*}
	Then define the map $\varphi: G \to G'$, where for all edges and halfedges $\varphi_E(x_i) \mapsto x$, $\varphi_H(x_i) \mapsto x$, and for the vertices $\varphi_V(x_i) \mapsto x$.  Then this is \'etale; but it is not an embedding, because $\varphi_V: V \to V'$ is not a monomorphism.
\end{eg}

\begin{defn}[Emb]
	Given a graph $G$, let $\widetilde{Emb(G)}$ be the collection of embeddings of graphs into $G$.  That is, $$\widetilde{Emb(G)} = \{f: H \to G \; |\; f \text{ is an embedding into }G\}$$.  
	
	Then $Emb(G)$ is the result of modding out by isomorphisms.  That is, if $f_1 \sim f_2$ then $[f_1]=\{f_1,f_2,\ldots\} \in Emb(G)$.
\end{defn}
Note that, rather than considering $Emb(G)$ to be a collection of maps, it can be considered the collection of graphs that can be embedded into $G$.
\begin{defn}[Vertex sum $\varsigma$]
	Let $f: G \to G'$ be an embedding in $Emb(G)$.  Then there is a corresponding element in the free commutative monoid on $V$, $$\sum_{v \in V}f(v) \in \N V'.$$  Denote this map by $\varsigma: Emb(G) \to \N V : f \mapsto \sum_{v \in V}f(v)$.
\end{defn}
In the above definition, note that the set of all $f(V)$ will be some subset of the vertices in $V'$; in fact, if $Y \subseteq X$, then $\sum_{v \in Y}v \leq \sum_{v \in X}v$.
 
\begin{figure}
	\label{eg:graphicalmapA}
	\centering
	\includegraphics[scale=0.8]{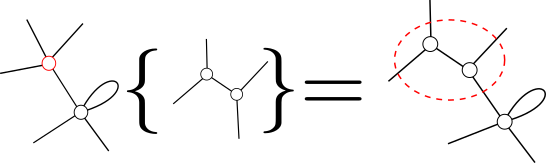}
	\caption{The graph in brackets in inserted into the red vertex to form the new graph.}
\end{figure}

Now that these definitions have been established, the definition of graphical map is presented.
\begin{defn}[Graphical map]
	\label{def:graphicalmap}
	A graphical map $\varphi : G \to G'$ consists of:
	\begin{itemize}
		\item A map of involutive sets $\varphi_0: E \to E'$, i.e.\ a map which respects the involutions, 
		\item A function $\varphi_1 : V \to Emb(G')$,
	\end{itemize}
	satisfying the following.
	\begin{enumerate}
		\item This inequality holds:
		$$\sum_{v \in V} \varsigma(\varphi_1(v)) \leq \sum_{w \in V'} w,$$
		\item There is a unique bijection making this diagram commute:
		\begin{center}
		\begin{tikzcd}
			nb(v) \arrow{r}{i} \arrow{d}{\cong} & E \arrow{d}{\varphi_0} \\
			\partial (\varphi_1(v)) \arrow{r}{inj} & E' \\
		\end{tikzcd}
		\end{center}
		Here the top map $i$ is the restriction of the involution
		and the bottom map is an injection. 
		Note that 
		$\partial (\varphi_1(v))$ refers to the boundary of the graph given by the image of the embedding $\varphi_1(v)$.
		\item If $\partial G = \emptyset$, then there exists a $v$ such that $\varphi_1(v)$ is not an edge. 
	\end{enumerate}
\end{defn}
The function $\varphi_1$ is essentially saying that for each vertex $v$ in $G$, a graph $H_v$ is assigned to it, which is then embedded into $G$ to get $G'$. 
In other words, one can select a vertex $v \in V(G)$ of degree $n$, and a graph $H$ with $n$ legs (that is, outer edges with only one edge connected to a vertex), and insert $H$ in place of $v$ to arrive at $G'$ (see Figure \ref{eg:graphicalmapA}).  This is written as $G\{H_v\}=G'$.  

The bijection is stating that the legs of $H_v$ correspond to the legs of the vertex $v$.
There is a slight overload of notation going on.  The expression $\partial(\varphi_1(v))$ actually refers to the boundary of the image of the embedding into $G'$, and the commutative diagram is saying that the boundary of the graph $H_v$ inserted into the vertex $v$ should match the boundary of $v$.

Then the inequality is stating that there should be fewer vertices in $G$ than in $G'$.

\begin{figure}
	\label{eg:maptoloop}
	\centering
	\includegraphics[scale=0.8]{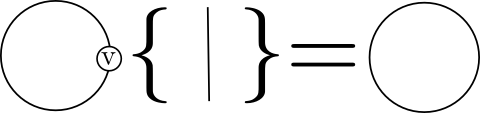}
	\caption{This map is not allowed by condition 3}
\end{figure}

The final requirement ensures that there is no map from the loop on one vertex to the graph consisting of a loop with no vertices  (see Figure \ref{eg:maptoloop}).  
Figure \ref{eg:graphicalmapA} shows an example of a graphical map.

		Graphical maps can be written down in multiple ways.  Firstly, in terms of two maps, one from edges to edges and the other from vertices to embeddings, as in Definition \ref{def:graphicalmap}.  Secondly, in terms of coface and codegeneracy maps, which shall be detailed in Appendix \ref{sec:facemaprelations}.  Thirdly, as graph substitutions.  A graph substitution, $G' = G\{H_v\}$, presents $G'$ as the result of replacing each vertex $v$ in $G$ with a graph $H_v$, where the boundary of $H_v$ is compatible with the neighbourhood of $v$.  According to the following lemma \citep[Prop 1.38]{robertson2019modular}, each graphical map can be written as a graph substitution followed by an inclusion.
		\begin{lem}
			\label{lem:graphicalmapgraphsub}
			Let $\varphi: G \to G'$ be a graphical map, and denote by $\varphi_v: H_v \inj G'$ each $\varphi_1(v) \in Emb(G')$.  Then there is an embedding $G\{H_v\} \inj G'$ which factors through each $\varphi_v$.
		\end{lem}

\begin{defn}[The category of graphs, $\Graph$]
\label{def:catGraph}
	Let $\Graph$ be the category with graphs as objects and graphical maps as morphisms, and with composition of graphical maps inherited from composition of graph substitutions. 
\end{defn}

Every graph provides a corresponding modular operad (see Definition \ref{def:GraphToModOp}), and therefore there is a functor from the category of graphs to the category of modular operads.  In the literature, $\Graph$ is used to refer to:
\begin{enumerate}
	\item The category of graphs.
	\item The functor from the category of graphs to the category of modular operads.
	\item The image of the functor.
\end{enumerate}
For clarity, in the following definition the functor will be referred to as $\Graph$ while the category of graphs will be denoted $\cat{Graph}$.  Afterwards, convention will be followed with $\Graph$ used for both.  Also, $G$ will refer to a graph while $\Graph(G)$ will refer to its corresponding modular operad.

\begin{defn}
	\label{def:GraphToModOp}
	There is a functor $\Graph: \cat{Graph} \to \cat{ModOpd}$ from the category of graphs to the category of modular operads.  Let $G$ be a graph in $\cat{Graph}$.  Then $\Graph(G)$ is defined as:
	\begin{itemize}
		\item Each arc of $G$ maps to a \hue{} in $ob(\Graph(G))$.  Then the \hueOrbits{} of $\Graph(G)$ are the edges of $G$.
		\item The operations are generated by the vertices.  Each vertex $v$ with boundary $\{c_1, \ldots, c_n\}$ generates a morphism for every permutation of boundary edges.  That is, for each $\sigma \in \Sigma_n$, $$v_{c_{\sigma(1)}, \ldots, c_{\sigma(n)}} \in \Graph(G)(c_1, \ldots, c_n).$$
			In addition, there is a unit $\eta_a \in \Graph(G)$ for each arc $a$.
		\item Two operations may be composed if they share at least one colour in common. 
		This corresponds in the graph to two subgraphs sharing an edge on the boundary of both.  If $w$ and $v$ are composed, then $w \circ v$ represents the subgraph of $G$ containing vertices $w$ and $v$, the edge between them, and any other arcs coming off them. 
		\item This composition is unital.  Each unit $\eta_e$ may be composed with any operation containing $\iota(e)$ as one of its \hues.  If $v$ is a subgraph containing the edge $\{e,i(e)\}$, then $v \circ \eta_e$ represents $v$ glued to $e$, which is thus left unchanged.  Similarly for $\eta_e \circ v$.
		\item The composition is associative.  The order in which edges in the graph are glued together does not matter. 
		\item A contraction may be performed on a vertex with two edges that are the same edge; that is, if it has a loop. 
	\end{itemize}
\end{defn}
\begin{eg}
	In the following graph, $G$, 
	\begin{center}
	\includegraphics[scale=0.8]{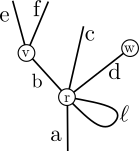}
	\end{center}
	the modular operad is defined as follows.
	\begin{itemize}
		\item The set of hues is $\{a_1, a_2,b_1,b_2,c_1,c_2,d_1,d_2,e_1,e_2,f_1,f_2,\ell_1,\ell_2\}$ (the set of colours is $\{a,b,c,d,e,f,\ell\}$).
		\item The operations are:
			\begin{align*}
				v_{e_2,f_2,b_1} & \in \Graph(G)(e_2,f_2,b_1) \\
				v_{e_2,b_1,f_2} & \in \Graph(G)(e_2,b_1,f_2) \\
				 & \vdots \\
				v_{b_1,e_2,f_2} & \in \Graph(G)(b_1,e_2,f_2) \\
				r_{a_1,b_2,c_2,d_2,\ell_2,\ell_1} & \in \Graph(G)(a_1,b_2,c_2,d_2,\ell_2,\ell_1) \\
				r_{a_1,b_2,c_2,\ell_1,d_2,\ell_2} & \in \Graph(G)(a_1,b_2,c_2,\ell_1,d_2,\ell_2) \\
				 & \vdots \\
				r_{d_2,\ell_1,c_2,\ell_2,b_2,a_1} & \in \Graph(G)(d_2,\ell_1,c_2,\ell_2,b_2,a_1) \\
				r_{a_1,b_2,c_2,d_2} & \in \Graph(G)(a_1,b_2,c_2,d_2) \\
				 & \vdots \\
				r_{d_2,c_2,b_2,a_1} & \in \Graph(G)(d_2,c_2,b_2,a_1) \\
				w_{d_1} & \in \Graph(G)(d_1) \\
			\end{align*}
			Along with all permutations of $v \circ_{(b_1)} r$, $r \circ_{(d_2)} w$, and $v \circ_{(b_1)} r \circ_{(d_2)} w$.
		\item There is a contraction, $\zeta (r_{a_1,b_2,c_2,d_2,\ell_2,\ell_1} ) = r_{a_1, b_2, c_2, d_2}$, along with all permutations thereof.
	\end{itemize}
\end{eg}

	\subsection{Graphical Sets}
			The aim is to define presheaves over the category of graphs, $\Graph$; infinity modular operads can be derived from there.  Graphical sets, also known as modular dendroidal sets, have been previously studied by \citet{robertson2019modular}, \citet{joyal2011feynman}, and \citet{sophiethesis, raynor2019distributivelaw}.

\begin{defn}[Graphical set]
	A graphical set $X$ is a presheaf $X: \Graph^{op} \to \cat{Set}$
\end{defn}
Denote by $X_G$ the set which is the image of a graph $G$ under $X$.
The category of graphical sets with natural transformations between them is denoted $\cat{gSet}$.
\begin{rmk}
In the same way that simplicial sets can be extended to simplicial objects, in the modular case there are also graphical objects, namely maps $\Graph^{op} \to \Cat$ where $\Cat$ is a cartesian monoidal category.  In particular, graphical groups will be used in this paper, presheaves $\Graph^{op} \to \cat{Grp}$.
\end{rmk}

The following definitions are examples of graphical sets that are particularly useful.
\begin{defn}[Representable]
	Given a graph $G \in \Graph$, the representable $\Graph[G] = \Graph(-,G)$.  That is, at each graph $H \in \Graph$, $\Graph[G]_H = \Graph(H,G)$.  
\end{defn}

Given a graph $G$, note the similarity in notation between its representable, $\Graph[G]$, and its associated modular operad, $\Graph(G)$ (see Definition \ref{def:GraphToModOp}).

For a full treatment of face maps, see Appendix \ref{sec:facemaprelations}.  In particular, note the definition of inner coface map of graphs, Definition \ref{def:graphicalInnerCoface}, which corresponds to a face map of modular operads.

\begingroup
\def\thetheorem{\ref{def:graphicalInnerCoface}}
\begin{defn*}[Inner coface map]
\sloppy
	Let $b$ be an inner edge of $G$.  Define $G/b$ as the graph created by contracting the edge $b$.  That is, if $x$ and $y$ are the vertices at either side of $b$, delete these and the edge $b$ and add a new vertex $xy$ with $nbhd(xy) = \lr{nbhd(x)\cup nbhd(y)} \setminus \{b\}$.  Then $\delta^b: G/b \to G$ is the associated inner coface map.

\begin{center}
	\includegraphics[scale=0.8]{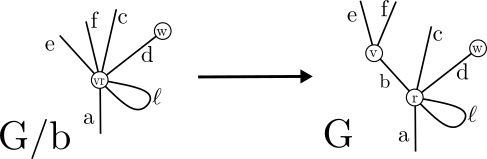}
\end{center}
\end{defn*}
\endgroup

In the appendix there are similar definitions for outer coface map  (Definition \ref{def:graphicalOuterCoface}) and the cosnip map (Definition \ref{def:graphicalCosnip}).  These coface maps (i.e. both inner and outer coface maps as well as snips), in the category of graphs, correspond to face maps when dealing with representables and modular operads.

\begin{defn}[Face]
	Let $G \in \Graph$ be a graph, and let $\alpha: H \to G$ be a coface map.  
	Then the $\alpha$-face of $\Graph[G]$ is the image of the map $\Graph[\alpha]: \Graph[H] \to \Graph[G]$.  It is denoted $\partial_\alpha \Graph[G]$.  
\end{defn}

Note that the \exedgename{} has no inner edges, and therefore has no inner faces.  Likewise, it has no vertices and therefore cannot have any outer faces, although it can itself be an outer face for any corolla.

\begin{defn}[Boundary of a representable]
	The graphical subset which is the union of all possible faces.
	$$\partial\Graph[G] = \bigcup_{\alpha} \partial_\alpha \Graph[G]$$
\end{defn}
\begin{defn}[Inner horn]
\label{def:InnerHorn}
	Let $\delta_e$ be an inner face map.  Then the inner horn $\Lambda_e$ is the graphical subset which is the union of all possible faces except $\delta_e$:
	$$\Lambda_e[G] = \bigcup_{\alpha \neq e} \partial_\alpha \Graph[G].$$
\end{defn}

There is a nerve defined in \citet{Robertson2019ModNerve}.
\begin{defn}[Modular nerve]
\label{def:ModularNerve}
	The nerve of a modular operad $\M$ is a functor $n_g: \cat{ModOpd} \to \cat{gSet}$.  This is given by, at each graph $G$,
	$$n_g(\M)_G = \cat{ModOpd}(\Graph(G), \M)$$
\end{defn}
Then, there is the modular realisation, a left adjoint to the nerve. 
\begin{defn}[Modular realisation]
	The aim is to define $\tau_g: \cat{gSet} \to \cat{ModOpd}$.
	Let $X$ be a graphical set.  Then a modular operad $\tau_g(X)$ can be defined from it as follows:
	\begin{itemize}
		\item The objects, $ob(\tau_g(X))$, are the set $X_\eta$
		\item The operations are generated by $\tau_g(X)(\underline{c})= X_{C_{\underline{c}}}$, where $C_{\underline{c}}$ is the corolla with arcs coloured by $\underline{c}$. 
		\item The face maps induce relations between these operations.
		\begin{enumerate}
			\item If $A$ is a \hue{} in $X_\eta$, and $\sigma: C_1 \to \eta$, then $\sigma^A(A) = id_A \in \tau_d(X)(A;A)$ is required, where $\sigma_A: X_{\eta} \to X_{C_1}$ is the map induced by $\sigma^A$.
			\item Let $T$ be the barbell; i.e. the tree with two vertices $v$ and $w$.   Then if $\delta_v$, $\delta_w$, and $\delta_{\overline{vw}}$ are similarly induced maps, require
		$$\delta_w(x) \circ_{\overline{vw}} \delta_v(x) = \delta_{\overline{vw}}(x)$$
		for all $x \in X_T$. 
		\end{enumerate}
	\end{itemize}
\end{defn}
Note that the barbell graph is somewhat special, as it is the simplest graph on two vertices, and its role in this definition is gluing two edges together. 
\begin{lem}
	The realisation $\tau_g$ is left adjoint to the nerve $n_g$.
\end{lem}
\begin{proof}
	The category $\Graph$ is small, and $\cat{ModOpd}$ has all small colimits.\citep{Robertson2019ModNerve}  
	Therefore, since the Yoneda embedding is full and faithful, the left adjoint is the left Kan extension of the inclusion $\Graph \inj \cat{ModOpd}$ along the Yoneda embedding \citep[Chapter X]{cwm}. 
\end{proof}

Note that the nerve of a graph does not always correspond to the representable of a graph (Remark 3.2 in \citep{Robertson2019ModNerve}), because the connection between morphisms in the category of graphs $\Graph$ and morphisms in the category of modular operads $\cat{ModOpd}$ is not straightforward.
As an example, consider the loop with one vertex, $K$.  This graph has two directed edges, $a$ and $b$, in opposing directions.  Consider the nerve of a graph $G$ at the point $K$.  Then
	$$ n_g(G)_K = \cat{ModOpd}(K,G)$$ and
	$$ \Graph[G]_K = \Graph(K,G)$$
But the former allows maps sending the directed edge $a$ to any directed edge $e \in G$ and the lone vertex of $K$ to the edge $e$, whereas the latter does not, so the result is $\Graph[G]_K \subsetneq n_g(G)_K$.

	\subsection{Segal Condition}
	\begin{defn}[Segal Core]
	\label{def:segalcore}
	Let $G$ be a graph with a least one vertex.  Then the Segal core $Sc[G]$ is the subset consisting of the union of corollae:
	$$Sc[G] = \bigcup_{v \in V(G)} \Graph [C_{\underline{c}(v)}]$$
	Where $\underline{c}(v)$ is the profile associated to the vertex $v$, and $C_{\underline{c}(v)}$ is the corolla associated to this vertex.  Sometimes the notation $C_v$ will be used to mean the same thing.
\end{defn}
Thus, one can write down a Segal condition that is similar to the alternative formulation of the Kan condition.  This is covered in \cite{Robertson2019ModNerve}.
\begin{defn}[Segal condition]
	A graphical set $X$ satisfies the Segal condition if, for all graphs $G \in \Graph$, there is a bijection:
	$$Hom(\Graph[G], X) \cong Hom(Sc[G],X)$$
	(induced by the Segal core inclusion $Sc[G] \inj \Graph[G]$).
\end{defn}
\section{The Inner Kan Condition}
    There are multiple ways of defining infinity categories, two of which are quasi-categories, which satisfy the Kan condition, and Segal categories.  The Segal condition has already been generalised to infinity modular operads (\citep{robertson2019modular,Robertson2019ModNerve}; now in this section the Kan condition is generalised, and a proof of equivalence is provided.
	\subsection{Quasi Modular Operads}
			Given its relevance in quasi operads and quasi categories, it is natural to examine what a Kan condition of modular operads should be.  This definition appears to be the correct one, since, as shown in Theorem \ref{thm:modularsegalkan}, it is equivalent to the Segal condition.  Note that the inner horn is defined in Definition \ref{def:InnerHorn}.
	\begin{defn}[Inner Kan condition]
		\label{def:kancond}
		A graphical set $X$ is said to satisfy the inner Kan condition if, for every graph $G$ and inner horn $\Lambda^e[G]$, the diagram
		\begin{center}
		\begin{tikzcd}
			\Lambda^e[G] \arrow{d}{} \arrow{r}{} & X \\
			\Graph[G] \arrow[dotted]{ur}{\exists!} & \\
		\end{tikzcd}
		\end{center}
		admits a filler.
	\end{defn}
	A quasi modular operad, therefore, can be defined in a similar way to a quasi category or a quasi operad.  
	\begin{defn}[Quasi modular operad]
		A quasi modular operad is a graphical set $X$ that satisfies the inner Kan condition.
	\end{defn}
	The reason for this definition becomes clear with Theorem \ref{thm:GraphicalNerve}, showing the link between quasi modular operads and modular operads.  If the filler in Definition \ref{def:kancond} is unique, then a graphical set is said to satisfy the \emph{strict} Kan condition, and has a corresponding modular operad.

	\subsection{Equivalence of Segal and Kan}
		\label{sec:SegalKanProof}
		The Segal condition for graphical sets (and its accompanying Nerve Theorem) was first introduced in \citet{Robertson2019ModNerve}.  
		In Theorem \ref{thm:modularsegalkan} I demonstrate an equivalence between the (strict) Kan condition and the Segal condition for modular operads, thus extending the nerve theorem to encompass the Kan condition.		
		
	Firstly, the simply connected case is proven in Example \ref{eg:segalkanbase}, Lemma \ref{lem:astersegalkan}, and Lemma \ref{thm:astersegalkanMainbody}.  Then Proposition \ref{prop:segalimplieskan} provides a proof that the Segal condition implies the Kan condition, and the combination of Proposition \ref{lem:spanningtree} and Theorem \ref{thm:modularsegalkan} completes the proof.  For an overview of face maps, see Appendix \ref{sec:facemaprelations}.
	
\renewcommand{\Aster}{\Graph}
\begin{eg}
	\label{eg:segalkanbase} 
	  The smallest tree containing at least three face maps, at least one being an inner face map, is the barbell, $B$ (see Definition \ref{def:barbell}).
	
	\begin{center}
		\includegraphics[scale=1]{./Images/BarbellEG}
	\end{center}
	
	The three face maps are $\delta_v$, $\delta_u$, and $\delta_e$, of which only $\delta_e$ is an inner face map.  Therefore, there is only one inner horn, 
	$$\Lambda^e[B] = \partial_u \Aster[B] \cup \partial_v \Aster[B].$$
	The result of $\partial_u$ is the corolla $C_v$, and likewise $\partial_v$ results in $C_u$.  Therefore, the inner horn is 
	$$\Lambda^e[B] = \Aster[C_v] \cup \Aster[C_u].$$
	This is the Segal core; therefore, in this instance, for any graphical set $X$, 
	$$\Hom(\Lambda^e[B],X) \cong  \Hom(Sc[B],X).$$
\end{eg}
Usually, the Segal core is merely a subset of the horn, and therefore $\Hom(\Lambda^e[G],X)$ will be a subset of $\Hom(Sc[G],X)$; in this case, however, they are in bijection.

\begin{lem}
	\label{lem:astersegalkan}
	Let $T$ be an unrooted tree with some inner edge $e$, and let $X$ be a graphical set.  If 
		$$\Hom(\Lambda^e[T],X) \cong \Hom(\Aster[T],X)$$
	then
		$$\Hom(Sc[T],X) \cong \Hom(\Aster[T],X).$$
\end{lem}
\begin{proof}
	This shall be proven via strong induction on the number of vertices in $T$.  Omitting those graphs with fewer than two vertices, the barbell graph shall form the base case, explored in Example \ref{eg:segalkanbase}.

	Given a map $f \in \Hom(Sc[T],X)$, we shall find a map $f_{\alpha} \in \Hom(Sc[T/\alpha],X)$ for each face in the horn.  Then, by the induction, we associate to it some function $h_\alpha \in \Hom(\Aster[T/\alpha], X)$.  Finally, note that $\Lambda^e[T] = \bigcup_{\alpha \neq e} \partial_{\alpha} \Aster[T]$, so $h \in \Hom(\Lambda^e[T],X)$ can be defined from the collection of maps $h_\alpha$, and then the Kan condition (Definition \ref{def:kancond}) gives a map $g \in \Hom(\Aster[T], X)$.

	Firstly, define each $f_\alpha$.  If $\alpha$ is an outer face map then $Sc[T/\alpha]$ includes into $Sc[T]$, since $Sc[T/\alpha]$ merely excludes the corolla associated to $\alpha$.  Otherwise, $\alpha$ is an inner face map (Defintion \ref{def:graphicalInnerCoface}.  Let $\overline{vw}$ denote the edge which is contracted, between the vertices $v$ and $w$, and let $vw$ denote the edge formed from the combination of $v$ and $w$.  For all other vertices, their associated corollae will be found in both $Sc[T]$ and $Sc[T/\alpha]$.  To complete the definition of $f_\alpha$, it is necessary to define $f_\alpha(\Aster[C_{vw}]$.  But by the inductive hypothesis $\Hom(\Aster[C_{vw}],X) \cong \Hom(Sc[C_{vw}],X)$, and $Sc[C_{vw}] = \Aster[C_v] \cup \Aster[C_w]$. 
	Let
	$$h(\partial_{\alpha} \Aster[T]) = h_{\alpha} (\Aster[T/\alpha]),$$
	noting that 
	we have $\partial_{\alpha} \Aster[T] = \Aster[T/\alpha]$.
	
	It remains to be shown that this is well defined, and that it is indeed a bijection. Let us consider well defined first.  Given two faces $\partial_{\alpha}$ and $\partial_{\beta}$, we require that
	$$h(\partial_{\alpha} \Aster[T]) = h(\partial_{\beta} \Aster[T])$$
	when restricted to the intersection $\partial_{\alpha} \Aster[T] \cap \partial_{\beta} \Aster[T]$.  But this follows from the identities (Appendix \ref{sec:facemaprelations}).  
	
	In most cases, $\delta^{\alpha} \delta^{\beta} = \delta^{\beta} \delta^{\alpha}$, so the aforementioned intersection is in bijection with the set $\partial_{\alpha} \partial_{\beta} \Aster[T]$.  
	However, if we consider an outer face map and an inner face map, corresponding to an inner edge incident to an outer vertex, then this conception of $\partial_{\alpha} \partial_{\beta} \Aster[T]$ does not entirely make sense.  Likewise with two outer face maps where one of them is done on a vertex that is not an outer vertex until the other face map is done.  In both of these cases, there is still an intersection involved, but we can't simply commute the two face maps.  Let our vertices be $p$ and $q$, connected via an edge $r$.  Then, due to the graphical identities, we know that $\partial_q \partial_p \Aster[T] = \partial_{pq} \partial_r \Aster[T]$, where $pq$ is the outer vertex resulting from contracting the edge $r$.  This intersection must be well defined too, because the graphical set satisfies these graphical identities.  (There are no faces $\partial_q$ and $\partial_{pq}$ of $\Aster[T]$, so we need not worry about them.) 
	
	Now to ensure that there is a bijection between the sets $\Hom(\Aster[T],X)$ and $\Hom(Sc[T],X)$.  Call the map defined above $$\varphi: \Hom(Sc[T],X) \to \Hom(\Aster[T],X).$$  Given any map in $\Hom(\Aster[T],X)$, it is clear that a map in $\Hom(Sc[T],X)$ can be defined via a restriction.  Call this map $$\rho: \Hom(\Aster[T],X) \to \Hom(Sc[T],X).$$  So it must be shown that $\varphi \circ \rho = id = \rho \circ \varphi$.  
	But, the graphical set relationship between $Sc[T]$, $\Lambda^e[T]$, and $\Aster[T]$ ensures that, as we go back and forth along $\rho$ and $\varphi$ between the homsets,
    we do indeed have $\varphi \circ \rho = id = \rho \circ \varphi$.

\end{proof}

\begin{lem}
	\label{thm:astersegalkanMainbody}
	Let $X$ be a graphical set.  If for all trees $T$ in $\Graph$ there is a bijection
	$$Hom(\Graph[T], X) \cong Hom(Sc[T],X)$$
	induced by the Segal core inclusion $Sc[T] \inj \Graph[T]$, then it is also true that, for every tree $T$ and inner horn $\lambda^e[T]$, the following diagram admits a unique filler:
		\begin{center}
		\begin{tikzcd}
			\Lambda^e[T] \arrow{d}{} \arrow{r}{} & X \\
			\Graph[T] \arrow[dotted]{ur}{\exists!} & \\
		\end{tikzcd}
		\end{center}
\end{lem}
\begin{proof}
	This can be shown directly.  Let $X$ be a graphical set that satisfies the condition above, $T$ be a tree with inner edge $e$, and let $h \in \Hom(\Lambda^e[T], X)$.  Each vertex of $T$ can be found in at least one face, so $h$ gives a map $f: Sc[T] \to X$.  And $f$ provides a map $g: \Aster[T] \to X$.  Therefore we have
	$$Hom(Sc[T],X) \subset Hom(\Lambda^e[T],X) \subset Hom(\Aster[T],X).$$
	
	We know that that $Hom(Sc[T],X) \cong Hom(\Aster[T],X)$ is a bijection, so we must therefore also have a bijection $\Hom(\Lambda^e[T],X) \cong \Hom(\Aster[T],X)$.
\end{proof}

\renewcommand{\Aster}{\Zhe}
The above lemma is analogous to Lemma \ref{thm:astersegalkan} in the Appendix, an equivalence of the Segal and Kan conditions for Astroidal sets.

Then, this theorem is extended to general graphs via the concept of spanning trees.   The core of this proof can be found in Lemma \ref{lem:spanningtree}, and it all comes together in Theorem \ref{thm:modularsegalkan}.

The definition of spanning tree used herein is slightly different to the standard definition in graph theory.
\begin{defn}[Spanning tree]
	Let $G$ be a graph.  Then a spanning tree of $G$, $T_G$, is a graph which
	\begin{enumerate}
		\item is a subgraph of $G$
		\item is connected
		\item has genus $0$
		\item includes each vertex of $G$, and the corolla associated to each vertex (that is, the set of edges connected to each vertex) is identical in both $G$ and $T_G$.
	\end{enumerate}
\end{defn}
Of particular note is that the standard definition of spanning graph does not have half edges, and so if the fourth property is true then all edges of the graph will be included and $T_G$ will be equal to $G$.  However, in this definition of graph half edges \emph{are} allowed, and so spanning trees correspond to breaking edges rather than removing them entirely.  Also, for any particular graph, note that there may be more than one spanning tree.  For example, consider this graph, $G$.
\begin{center}
	\includegraphics[scale=0.5]{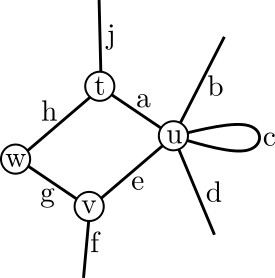}
\end{center}
It has four different spanning trees, each corresponding to breaking one edge in the cycle (in addition to breaking the loop).  Two of them are:
\begin{center}
	\includegraphics[scale=0.4]{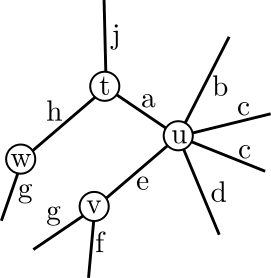} \hspace{1cm}
	\includegraphics[scale=0.4]{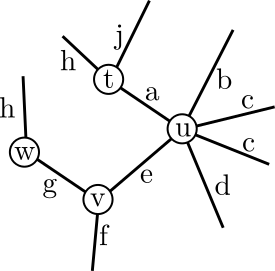}
\end{center}
Note that, in each case, the loop (edge $c$) must also be broken.  Additionally, only one edge in the cycle may be removed, or the resulting tree will not be connected.  For each broken edge, both sides have been labelled by the label of the original edge; however, it is important to note that each edge consists of two half edges, and the two halves will each be one of these half edges.

Also, note that the Segal cores of a graph and its spanning trees will be identical.
\begin{prop}
	\label{lem:spanningtree} 
	Let $X$ be a graphical set, and consider a graph $G$.  Let $T_i$ be a collection of spanning trees of $G$ such that each edge in $G$ is included in at least one spanning tree.  
	If each of these spanning trees satisfy $$\Hom(Sc[T_i], X) \cong \Hom(\Graph [T_i], X),$$ then 
	$$\Hom(Sc[G], X) \cong \Hom(\Graph [G], X)$$
\end{prop}
\begin{proof}
	We wish to construct a bijection $\Hom(Sc[G], X) \cong \Hom(\Graph [G], X)$.  Given a map $f: Sc[G] \to X$, we wish to find a map $h: \Graph [G] \to X$, then show that $f \mapsto h$ is a bijection.
	
	Firstly, consider any particular spanning tree $T_i$.  Since $Sc[G] \cong Sc[T_i]$, let $f_i: Sc[T_i] \to X$ be defined by $f_i(x) = f(x)$.  Then, by our hypothesis, we can construct a map $h_i: \Graph[T_i] \to X$.  
	Then, each $T_i$ can be reached from $G$ by a series of cosnip maps.  We also know that the cosnip is the only face or degeneracy map which changes genus.   So we get $\Graph [G] \to \Graph [T_i]$ via a series of snips.  These shall be denoted $d_i$ in this proof.
	
	From this collection of maps $h_i$ and $d_i$, $h$ shall be constructed as $h = \bigcup h_i \circ d_i$.  It remains to be shown that this map is well defined, and a bijection.  Let this map $f \mapsto h$ be denoted by $\psi$.  	The following diagram is illustrative.
	\begin{center}
	\begin{tikzcd}
		\Graph[G] \arrow{d}{d_i} \arrow{dr}{h} & \\
		\Graph[T_i] \arrow{r}{h_i} \arrow{d}{\cong} & X \arrow{d}{=}\\
		Sc[T_i] \arrow{d}{=} \arrow{r}{f_i} & X \\
		Sc[G] \arrow{ur}{f} & \\
	\end{tikzcd}
	\end{center}
	
	This map is indeed well-defined.  Consider two spanning trees $T_j$ and $T_k$. We wish to compare $h = h_j \circ d_j$ and $h' = h_k \circ d_k$.  But the coface and codegeneracy maps in $\Graph$  are associated with the face and degeneracy maps in $X$ according to the graphical relations. That is, if $x \in \Graph[G]_S$, then $x: S \to G$ maps to some $d_j(x)$, which is $x$ restricted to the spanning tree $T_j$.  Because each map $\Graph[T_i] \to X$ is induced by the coface maps $d^i$,  $h$ and $h'$ agree on vertices and their neighbouring half-edges.  Then, the graphical relations ensure that these maps agree on edges between two different corollae. 
	
	Now, consider a map $h \in \Hom(\Graph[G],X)$.  From $h$, a map $f \in \Hom(Sc[G],X)$ can be defined by restricting $h$ to each corolla $\Graph[C_v]$.  This provides a map $\psi^{-1}$.  

	Then $\psi \circ \psi^{-1} = \psi^{-1} \circ \psi = id$ because the information at each corolla is preserved.  Thus, the map $\psi: \Hom(Sc[G],X) \to \Hom(\Graph[G],X)$ is a bijection.
\end{proof}

\begin{prop}
	\label{prop:segalimplieskan}
	Let $X$ be a graphical set that satisfies the graphical Segal condition.  Then it also satisfies the strict inner graphical Kan condition.
\end{prop}
\begin{proof}
	Let $X$ be a graphical set that satisfies the Segal condition.  
	
	Firstly, consider the following series of inclusions (for all graphs $G$ and all inner faces $e$)
	$$Sc[G] \subset \Lambda^e[G] \subset \Graph[G]$$
	An element of $Sc[G]_H$ is an element of $\Graph[C_v]_H$, where $C_v$ is the corolla associated to a vertex $v$ of $G$.  That is, a graphical map $H \to C_v$.  But such a graphical map partially defines a graphical map $H \to \partial_\alpha \Graph[G]$, because each corolla is contained in at least one face of $G$.  Therefore, $Sc[G] \subset \Lambda^e[G]$.  Likewise, each graphical map in the inner horn partially defines a graphical map on the whole graph, so $\Lambda^e[G] \subset \Graph[G]$.
	
	Therefore, if we have a map $Sc[G]\to X$, and thus by the Segal condition one $\Graph[G] \to X$, we must also have one $\Lambda^e[G] \to X$.
	Therefore, we have 
	$$Hom(Sc[G],X) \subset Hom(\Lambda^e[G],X) \subset Hom(\Graph[G],X).$$
	
	The Segal condition states that $Hom(Sc[G],X) \cong Hom(\Graph[G],X)$ is a bijection, so we must therefore also have a bijection $\Lambda^e[G] \cong \Graph[G]$, giving us the Kan condition.  
\end{proof}

\begin{thm}
	\label{thm:modularsegalkan}
	Let $X$ be a graphical set.  Then $X$ satisfies the graphical Segal condition if and only if $X$ satisfies the strict inner graphical Kan condition.
\end{thm}
\begin{proof}
	The proof that Segal implies Kan is found in Proposition \ref{prop:segalimplieskan}.
	Now assume that $X$ satisfies the Kan condition.  We shall prove that 
	\begin{align*}
		\Hom(\Lambda^e[G],X) & \cong \Hom(\Graph[G],X) \\
		\implies \Hom(Sc[G],X) & \cong \Hom(\Graph[G],X) \\
	\end{align*}
	for all graphs $G$.  
	
	Assume $G$ is a higher genus graph, and consider a spanning tree $T_G$.  Since $\Hom(\Lambda^e[G],X) \cong \Hom(\Graph[G],X)$, we know that $\Hom(\Lambda^e[T_G],X) \cong \Hom(\Graph[T_G],X)$ for all spanning trees $T_G$, unless $T_G$ would remove the edge $e$.  In this case, $\Lambda^e[T_G],X)$ is undefined.  However, a spanning tree may only remove the inner edge $e$ if it exists as part of some cycle, but that implies that there must be some other edge in the cycle, $f$, and the spanning tree(s) that remove $f$ will include $e$.  
	Note that, if $e$ is part of a cycle, then there may be a spanning tree $T_G$ which is formed by removing $e$.  In this case, $\Lambda^e[T_G]$ cannot be defined; however, since $e$ is part of a cycle there exists some other edge $e'$.
	By Lemma \ref{thm:astersegalkanMainbody}, we know that, for each spanning tree, $\Hom(Sc[T_G],X) \cong \Hom(\Graph[T_G],X)$, and thus by Lemma \ref{lem:spanningtree} we have $\Hom(Sc[G],X) \cong \Hom(\Graph[G],X)$, as required.
\end{proof}
\subsection{The Nerve Theorem}
Being able to move between the graphical sets and modular operads is important, since the definition of infinity modular operads is usually given in terms of graphical sets, just as infinity categories are defined in terms of Simplicial sets.  Quasi categories were first studied by \citep{boardman2006homotopy}, and the equivalence between the Segal and Kan conditions is found in \citep{grothendieck1960technique} and \citep{segal1968classifying}.  This was further expanded to Dendroidal sets (Corollary 2.6 in \cite{cisinski2013dendroidal}) and Astroidal sets (Theorem 6.7 in \cite{hackney2019higher}, Theorem 4.3.7 in \cite{mythesis}), and I now extend it to graphical sets.
	\begin{thm}[Graphical Nerve Theorem]
	\label{thm:GraphicalNerve}
    \sloppy
	Let $X$ be a graphical set ${X : \Graph^{op} \to \cat{Set}}$.  Then the following are equivalent:
	\begin{enumerate}
		\item There exists a modular operad $\M$ such that $n_g(\M) = X$
		\item $X$ satisfies the (strict) inner Kan condition
		\item $X$ satisfies the Segal condition
	\end{enumerate}
\end{thm}
The equivalence between 1 and 3 can be found in \citep[Theorem 3.6]{Robertson2019ModNerve}.  The proof of the equivalence between 2 and 3 is original and is found above in Theorem \ref{thm:modularsegalkan}.

\newpage

\begin{appendices}

\section{Cyclic Operads}
\label{app:cycOp}
Many of the results in this paper also apply to cyclic operads, as they are genus $0$ modular operads.

	In \citep{gretzler1994cyclic}, Getzler and Kapranov introduce the notion of a cyclic operad, as a generalisation of the cyclic homology of associative algebras.  At the time, cyclic operads were important in the study of homotopy theory \citep{may2006geometry} and topological field theory \citep{getzler1994batalin}.  They are also of interest for their connections to physics and dagger categories \citep{baez2006quantum} and the Grothendieck-Teichm\"uller group \citep{de2017operads}.  Cyclic operads are also a generalisation of dagger categories, in the same way that operads are a generalisation of categories.

Infinity cyclic operads were first defined in \citep{hackney2019higher} using a Segal condition, and a model structure for simplicially enriched cyclic operads is detailed in \citep[Section 6]{drummondcole2018dwyerkan}).  This paper provides a framework for quasi cyclic operads.  For a detailed picture of quasi cyclic operads, see \citet{mythesis}.

\subsection{Cyclic Operads}

The first definition of a cyclic operad is due to \citet{gretzler1994cyclic}, with a new axiom added by Van der Laan \citep[Section 11]{van2004operads}.  There is a definition by \citet{obradovic2016monoid} based on a monoidal construction, while  Hackney, Robertson, and Yau give a definition based on monads in Section 5.1 of \citep{hackney2019higher}.  Definition \ref{def:cyclicOperad} is equivalent to this last one. 

Note that some definitions of cyclic operads 
involve dualising the objects \citep{cheng2014cyclic, drummondcole2018dwyerkan}, but that is not used in this paper.

\begin{defn}[Cyclic operad] 
\label{def:cyclicOperad}
A (coloured, symmetric) cyclic operad $\Cat$ over a closed symmetric monoidal category $\mathcal{E}$ is defined with the following properties.
	\begin{itemize}
		\item A set of objects, $ob(\Cat)$, sometimes called colours.
		\item For each profile over $ob(\Cat)$, $\underline{c} = (c_1, \ldots, c_n)$, there is an object of $\mathcal{E}$ $\Cat(c_1, \ldots c_n)$.  If $\M$ is $\cat{Set}$, then this is the set of operations.
		\item There is a right action of the symmetric group; for any profile $\underline{c}$ and permutation $\sigma \in \Sigma_n$, there is a bijection $\Cat(c_1, \ldots c_n) \to \Cat(c_{\sigma(1)}, \ldots c_{\sigma(n)})$.
		\item There is a composition operation.  Given any two operations $\alpha \in \Cat(\underline{c})$ and $\beta \in \Cat(\underline{d})$ where the final object in $\underline{c}$ is equal to the $i$th object in $\underline{d}$, there is a composition $\alpha \circ_i \beta$  which is a map
			$$\Cat(\underline{c}) \otimes \Cat(\underline{d}) \to  \Cat(\underline{cd}_i)$$ 
		Where, for the sake of brevity, $$\underline{cd}_i = (d_1, \ldots d_{i-1}, c_1, \ldots c_{m-1}, d_{i+1}, \ldots d_n).$$
		\begin{itemize}
		\item The composition is associative.  That is, all such diagrams commute:
		\begin{center}
			\begin{tikzcd}
				\Cat(\underline{a})\otimes \Cat(\underline{b})\otimes \Cat(\underline{c}) \arrow{r}{id \otimes \circ_j} \arrow{d}{\circ_i \otimes id} & \Cat(\underline{a}) \otimes \Cat(\underline{bc}_j)\arrow{d}{\circ_{i+j-1}} \\
				\Cat(\underline{ab}_i) \otimes \Cat(\underline{c}) \arrow{r}{\circ_j} & \Cat(\underline{abc}_{ij})
			\end{tikzcd}
		\end{center}
		where the lengths of $\underline{a}$, $\underline{b}$, and $\underline{c}$ are $\ell$, $m$, and $n$, respectively.
		\item Composition is unital.  That is, for each colour $c$ there exists an identity $\eta_c$, such that for all $\theta \in \Cat(\underline{c})$, $\theta \circ_1 \eta_{c_n} = \theta = \eta_{c_i} \circ_i \theta$.
		\item Composition is equivariant.  That is, it commutes with the action of the symmetric group.  So for any two morphisms $\alpha \in \Cat(\underline{c})$ and $\beta \in \Cat(\underline{d})$, and any $\sigma, \tau \in \Sigma_n$, 
		$$\alpha \circ_{\sigma(i)} \sigma(\beta) = \sigma'(\alpha \circ_i \beta )$$
		$$\tau(\alpha) \circ_i \beta = \tau'( \alpha \circ_i \beta)$$
		Where $\sigma' \in \Sigma_{n+m}$ refers to the element that acts on $$\{ d_1, \ldots, d_{i-1}, c_1, \ldots, c_m, d_{i+1}, \ldots, d_n \}$$ by doing $\sigma$ on each $d_k$ and the identity on each $c_j$, and $\tau' \in \Sigma_{n+m}$ does the identity on each $d_k$ and permutes each $c_j$ according to $\tau$.
		\end{itemize}
	\end{itemize}
	Then $\Cat$ is a (coloured symmetric) cyclic operad.
\end{defn}

\begin{defn}
	Let $\Op$ and $\oP$ be two cyclic operads.  A morphism $f: \Op \to \oP$ consists of the following information:
	\begin{itemize}
		\item A function $f': ob(\Op) \to ob(\oP)$ that sends the colours of one operad to the set of colours of the other.
		\item For each profile $\underline{c}$, a map $f_{\underline{c}}: \Op(\underline{c}) \to \oP(f'(\underline{c}))$ that commutes with the structure maps.  i.e.
		\begin{itemize}
			\item Identity: $f(id_\eta) = id_{f(\eta)}$
			\item Symmetric group action: $f(\sigma(v)) = \sigma(f(v))$ 
			\item Composition: $f_{\underline{cd}_i}(v \circ_i u) = f_{\underline{c}}(v) \circ_i f_{\underline{d}}(u)$
		\end{itemize}
		Here $v \in \Op(\underline{c})$ and $u \in \Op(\underline{d})$
	\end{itemize}
\end{defn}
Then the category $\cat{CycOpd}$ can be defined as the category with cyclic operads as objects and morphisms as defined above.

\subsection{Cyclic Dendroidal Sets}
\subsubsection{The Category of Unrooted Trees}
\begin{defn}[Unrooted tree]
	\label{def:unrootedtree}
	Let $G$ be a simply connected graph.  Let $I$ be a subset of the set of outermost vertices, where the outermost vertices are those connected to only one edge.  Then $(G,I)$ is an unrooted tree, usually shortened to just $G$.  As a convention, the input vertices are not drawn, and the word ``vertex'' refers only to the remaining vertices. 
\end{defn}

Morphisms are compositions of coface and codegeneracy maps.  
There are two types of coface maps, inner and outer.  In the following, let $T$ be an unrooted tree. 
\begin{defn}[Inner coface map]
	Let $c$ be an inner edge of $T$.  Define $T/c$ as the tree created by contracting the edge $c$.  That is, if $x$ and $y$ are the vertices at either side of $c$, delete these and the edge $c$ and add a new vertex $xy$ with $nbhd(xy) = \lr{nbhd(x)\cup nbhd(y)} \setminus \{c\}$.  Then $\delta^c: T/c \to T$ is the associated inner coface map.
\end{defn}

\begin{defn}[Outer coface map]
	Let $u$ be an outer vertex of $T$.  Then define $T/u$ to be the tree with $u$ and all its leaves deleted.  Then $\delta^u: T/u \to T$ is the associated outer coface map.
\end{defn}

\begin{defn}[Codegeneracy map]
	Let $c$ be any edge of $T$.  Then define $T_c$ to be the graph where the edge $c$ has been subdivided with exactly one vertex.  Then $\delta^c: T_c \to T$ is the associated codegeneracy map.
\end{defn}

\begin{eg}
	The single edge $\eta$ can be included into an $n$-leaved corolla $C_n$ in $(n+1)$ ways.
\end{eg}

These satisfy relations similar to those in the categories of dendrices and simplices \citep{mythesis}.

\subsubsection{Astroidal sets}
An astroidal set is a collection of sets patterned after the category of astrices. 
\begin{defn}[Astroidal Set]
	An astroidal set is a presheaf $X: \Aster^{op} \to \cat{Set}$.
\end{defn}
In other words, an astroidal set $X$ consists of the following information.
\begin{itemize}
	\item For each $T \in \Aster$, a set $X(T)$, denoted $X_T$.  Note that $X_T$ is called the set of astrices of shape $T$,
	\item For each morphism $f: S \to T$, a function $X_f: X_T \to X_S$,
	\item With $X_{id: T \to T} = id: X_T \to X_T$,
	\item Given two morphisms
		\begin{tikzcd}
			R \arrow{r}{\beta} & S \arrow{r}{\alpha} & T 
		\end{tikzcd}
		in $\Aster$, $X_{(\alpha \circ \beta)} = X_\beta \circ X_\alpha$.
\end{itemize}

\begin{defn}[Morphisms]
	Let $X: \Aster^{op} \to \cat{Set}$ and and $Y: \Aster^{op} \to \cat{Set}$ be two astroidal sets.  Then a map between them $f: X \to Y$ is a natural transformation.  In other words,
	\begin{itemize}
		\item For each tree $T \in \Aster$, there is a map $f_T : X_T \to Y_T$,
		\item If $\alpha: S \to T$ is a morphism in $\Aster$, then the following diagram commutes.
		\begin{center}
			\begin{tikzcd}
				X_t \arrow{r}{X_{\alpha}} \arrow{d}{f_T} & X_S \arrow{d}{f_S} \\
				Y_T \arrow{r}{Y_{\alpha}} & Y_s
			\end{tikzcd}
		\end{center}
	\end{itemize}
\end{defn}
The category of astroidal sets is defined in the expected way.
\begin{defn}[Astroidal sets]
	The category of astroidal sets, $\cat{aSet}$ is defined to be the category consisting of $ob(\cat{aSet}) = \{ \Aster^{op} \to \cat{Set} \}$ with natural transformations as morphisms. 
\end{defn}

\subsection{Quasi cyclic operads}
\subsubsection{The Kan condition}
\begin{defn}[Inner Kan complex]
	\label{def:asterKan}
	Let $X$ be a cyclic dendroidal set, $f: \Lambda^k[T] \to X$ be an inner horn, and let $j: \Lambda^k[T] \rightarrowtail \Aster[T]$ be the inclusion.  Then a filler for $f$ is a map $g : \Aster[n] \to X$ such that $f = g \circ j$.
	\begin{center}
	\begin{tikzcd}
		\Lambda^k[T] \arrow{r}{f} \arrow{d}{j} & X \\
		\Aster[T] \arrow[dotted]{ur}{g} & \\
	\end{tikzcd}
	\end{center}
	The cyclic dendroidal set $X$ is said to be an inner Kan complex if every inner horn has a filler.
\end{defn}

Then, a cyclic dendroidal set satisfying the Kan condition is an inner Kan complex, otherwise known as a quasi cyclic operad.
If the filler is unique 
then it is a strict inner Kan complex;
that is, there is a corresponding cyclic operad (see Theorem \ref{thm:asterNerve}).

\subsubsection{The Segal condition}

The Segal condition for cyclic dendroidal sets shall be defined analogously to the Segal condition for dendroidal sets.  A simplicially enriched version is studied in \citet[Definitions 1.34 and 8.8]{hackney2019higher}. 

\begin{defn}[Segal core]
	\label{def:segalcoreast}
	Let $T$ be a tree with at least one vertex.  Recall the corolla $C_n$.  The Segal core $Sc[T]$ is the sub-object of $\Aster [T]$ defined as the union of all the images of maps $\Aster[C_n] \to \Aster[T]$ corresponding to sub-trees of shape $C_n \to T$.  
Note that such a map is completely determined, up to isomorphism, by the vertex $v$ of $T$ in its image.  Let $n(v)$ be the number of input edges that $v$ has.  Therefore, one can write $$Sc[T] = \bigcup_{v \in V(T)} \Aster [C_{n(v)}].$$
\end{defn}

Then the Segal condition can be defined as follows.
\begin{defn}
	A cyclic dendroidal set $X$ satisfies the Segal condition if for every tree $T$ the map 
	$$\Hom(Sc[T], X) \to \Hom(\Aster[T],X)$$
	is a bijection.
\end{defn}

As an example, which will serve to partially illustrate the following theorem, consider the barbell graph (Definition \ref{def:barbell}).
\begin{eg}
	\label{eg:graphicalsegalmaps}
	Let $B$ be the tree containing two vertices, $v_1$ and $v_2$, joined by an edge $e$, where $v_1$ may have a different number of legs than $v_2$.  Consider the astroidal set $\Aster[B]$.  This satisfies the Segal condition as written above.  It also satisfies an alternative characterisation of the Segal condition: for every tree $T$ there is a bijection
		$$ \Aster[B]_T \cong \Aster[B]_{C_{v_1}} \times_{\Aster[B]_\eta} \Aster[B]_{C_{v_2}}.$$
\end{eg}
\begin{proof}
	The Segal map shall be examined for the graph $B$ in relation to this graphical set.  For most graphs $H$ unrelated to $B$, the set $\Aster[B]_H$ will be empty (the main exceptions being degeneracy maps), and thus uninteresting for the purpose of example.
	
	Consider the following square.
	\begin{center}
	\begin{tikzcd}
		G & C_{v_1} \arrow{l}{} \\
		C_{v_2} \arrow{u}{} & \eta \arrow{u}{} \arrow{l}{} \\
	\end{tikzcd}
	\end{center}

	For any tree $H$, there are many choices for the map $\eta \inj H$, one for each edge in $H$.  However, in the above diagram the only way for it to commute is to chose the maps $\eta \inj C_{v_1}$ and $\eta \inj C_{v_2}$ where $\eta$ is mapped to their edge in common, $e$.
	
	Consider the components of the Segal condition:
	\begin{align*}
		\Aster[B]_\eta & = \Hom(\eta, B) \\
		\Aster[B]_{C_{v_1}} & = \Hom(C_{v_1}, B) \\
		\Aster[B]_{C_{v_2}} & = \Hom(C_{v_2}, B) \\
	\end{align*}
	The first set has already been discussed.  For the other two, there is only one option, the outer coface maps $\delta^{v_2} : C_{v_1} \to B$ and $\delta^{v_1} : C_{v_2} \to B$.
	
	Although in theory there may be multiple choices of map $\Aster[B]_{C_{v_1}} \to \Aster[B]_\eta$ and $\Aster[B]_{C_{v_2}} \to \Aster[B]_\eta$, the way that $B$ is connected means that only the map induced by the above commutative diagram can be used.
	
	Thus, the pullback
	\begin{center}
	\begin{tikzcd}
		\Aster[B]_{C_{v_1}} \times_{\Aster[B]_\eta} \Aster[B]_{C_{v_2}} \arrow{r}{} \arrow{d}{} & \Aster[B]_{C_{v_1}} \arrow{d}{} \\
		\Aster[B]_{C_{v_2}} \arrow{r}{} & \Aster[B]_\eta \\
	\end{tikzcd}
	\end{center}
	contains only one map, the map which sends $C_{v_1} \to C_{v_1}$ and $C_{v_2} \to C_{v_2}$, up to their equivalence on the edge $e$.  But this describes the only map in $\Aster[B]_B$, $B \to B$.

\end{proof}

Just as astroidal sets can be defined as functors $\Aster^{op} \to \cat{Set}$, astroidal groups can be defined as functors $\Aster^{op} \to \cat{Grp}$.  The following theorem is proven here for astroidal groups, but is conjectured to also hold for graphical groups. 
\begin{prop}
	\label{thm:GroupsAreSegal}
	Let $X$ be an astroidal group.  If, for all trees $T$, 
	$$\bigcap_{v \in V(T)} \ker( \xi_{v} ) = \{1\},$$
	then the underlying astroidal set satisfies the Segal condition.  Note that each map $\xi_{v}: \Aster[T] \to \Aster[C_v]$ refers to a series of outer face maps resulting in just the corolla remaining.
\end{prop}
\begin{proof}
	Firstly, this is true for the barbell, by Example \ref{eg:graphicalsegalmaps}.  The proof will proceed by induction, adding in an extra outer vertex to the tree.  Let $T_{n-1}$ denote the tree $T$ without the  outer vertex $v_n$, so that $T$ has $n$ vertices and $T_{n-1}$ has $n-1$ vertices.  Assume the theorem holds for all trees up to $n-1$ vertices.  
	
	Consider the following diagram, where the maps $X_{C_{v_i}} \to X_\eta$ are induced by $\eta \to C_{v_i}$, and $Sc_{n-1}$ is the previous iteration of the map $Sc_n$.
	\begin{center}
	\begin{tikzcd}
		& & X_{T_n} \arrow{dr}{} & \\
		X_T \arrow[bend left=20]{urr}{Sc_{n-1}} \arrow[bend right=20]{drr}{\xi_{v_n}} \arrow{r}{Sc_n} & X_{T_{n-1}} \times_{X_\eta} X_{C_{v_n}} \arrow{ur}{} \arrow{dr}{} & & X_\eta \\
		& & X_{C_{v_{n}}} \arrow{ur}{} & \\
	\end{tikzcd}
	\end{center}

	We wish to show that $Sc_n$ is a bijection.  
	This map is obviously surjective, since $X_{T_{n-1}} \times_{X_\eta} X_{C_{v_n}}$ is a subset of $X_T$.
	
	Now, in order to show injectivity, assume for the sake of contradiction that there is some element $a \in ker(Sc)$ which is not the identity.  It is known that $X_{T_{n-1}} \times_{X_\eta} X_{C_{v_n}}$ is the limit
	so we know that for all $i$,  $\xi_{v_n}(a) = id$.  Therefore $a \in \bigcap_{v \in V(T)} \ker( \xi_{v})$.  But our hypothesis states that the identity is the only element in this set, a contradiction.  Therefore our function is injective.
	
	Now to show that this satisfies the Segal condition.  Consider the following three equivalences (given by the Yoneda Lemma) 
	\begin{align*}
		X_T & \cong \Hom(\Aster[T],X) \\
		X_{C_v} & \cong \Hom(\Aster[C_v],X) \\
		X_{\eta} & \cong \Hom(\Aster[\eta],X) \\
	\end{align*}
	as well as the pullback $X_{T_{n-1}} \times_{X_\eta} X_{C_{v_n}}$.  It suffices to show that $X_{T_{n-1}} \times_{X_\eta} X_{C_{v_n}} \cong \Hom(Sc[G],X)$.
	
	By Definition \ref{def:segalcore}, the Segal core is the union of the images of the maps $\Aster[C_n] \to \Aster[T]$ corresponding to sub-trees of shape $C_v \to T$.  But these maps will agree whenever two corollae share an edge.  Likewise, the pullback $X_{T_{n-1}} \times_{X_\eta} X_{C_{v_n}}$ corresponds to maps $X_{C_v} \cong \Hom(\Aster[C_v],X)$ which agree over $X_{\eta}$.  
\end{proof}

\subsubsection{The Nerve Theorem}
\begin{lem}
	\label{thm:astersegalkan}
	Let $X$ be an astroidal set.  If $X$ satisfies the astroidal Segal condition then it satisfies the strict inner astroidal Kan condition.
\end{lem}
\begin{proof}
	This can be shown directly.  Let $X$ be an astroidal set that satisfies the Segal condition, $T$ be a tree with inner edge $e$, and let $h \in \Hom(\Lambda^e[T], X)$.  Each vertex of $T$ can be found in at least one face, so $h$ gives a map $f: Sc[T] \to X$.  By the Segal condition, $f$ provides a map $g: \Aster[T] \to X$.  Therefore we have
	$$Hom(Sc[T],X) \subset Hom(\Lambda^e[T],X) \subset Hom(\Aster[T],X).$$
	
	The Segal condition states that $Hom(Sc[T],X) \cong Hom(\Aster[T],X)$ is a bijection, so we must therefore also have a bijection $\Hom(\Lambda^e[T],X) \cong \Hom(\Aster[T],X)$, giving us the Kan condition
\end{proof}

\begin{thm}[Nerve theorem]
	\label{thm:asterNerve}
	Let $\mathcal{A}$ be an astroidal set.  Then the following are equivalent:
	\begin{enumerate}
		\item There exists a cyclic operad $\Cat$ such that $\mathcal{A} \cong n_a(\Cat)$
		\item $\mathcal{A}$ satisfies the Segal condition
		\item $\mathcal{A}$ satisfies the strict inner Kan condition
	\end{enumerate}
\end{thm}
This is a corollary of Theorem \ref{thm:GraphicalNerve}.  
The proof of the equivalence of $1$ and $2$ can be found in \citet[Theorem 6.7]{hackney2019higher}, and Lemma \ref{thm:astersegalkan} shows that $2$ implies $3$.  
	
\section{Face and Degeneracy Maps in the Category of Graphs}
	\label{sec:facemaprelations}
	
In the category of simplicies $\Delta$, along with other related categories such as $\Omega$ for trees, it was possible to decompose the morphisms into faces and degeneracies.  The morphisms in the category of graphs $\Graph$ can also be generated by some elementary morphisms, which can be arranged into face maps and then degeneracy maps (Lemma \ref{lem:graphtofacedegen}).  First, these face and degeneracy maps must be defined.  There are two types of coface maps, as well as degeneracy maps.  In the following, let $G$ be a graph.  

\begin{defn}[Inner coface map]
\label{def:graphicalInnerCoface}
	Let $b$ be an inner edge of $G$.  Define $G/b$ as the graph created by contracting the edge $b$.  That is, if $x$ and $y$ are the vertices at either side of $b$, delete these and the edge $b$ and add a new vertex $xy$ with $nbhd(xy) = \lr{nbhd(x)\cup nbhd(y)} \setminus \{b\}$.  Then $\delta^b: G/b \to G$ is the associated inner coface map.
\end{defn}
\begin{center}
	\includegraphics[scale=0.8]{./Images/InnerFaceEg}
\end{center}
And the same example as a graph inclusion, where a barbell graph is inserted into the appropriate vertex.
\begin{center}
	\includegraphics[scale=0.8]{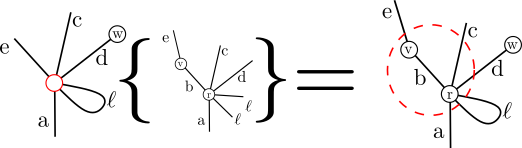}
\end{center}
That is, $\delta^b: G/b \to G$ is defined with $\delta^b_0$ the identity on all edges in $G/b$ and $\delta^b_1$ given by 
\begin{align*}
	w & \mapsto C_w \\
	vr & \mapsto B \\
\end{align*} 
where $B$ is the barbell graph given in the braces in the image.  The two properties of the graphical maps definition hold because the total vertices in the embeddings matches the total vertices in $G$, and each border matches where it attaches to its neighbouring vertices.

\begin{defn}[Outer coface map]
\label{def:graphicalOuterCoface}
	Let $v$ be an outer vertex of $G$. 
	Then define $G/v$ to be the graph with $v$ and all its leaves deleted.  Then $\delta^v: G/v \to G$ is the associated outer coface map.
\end{defn}
\begin{center}
	\includegraphics[scale=0.8]{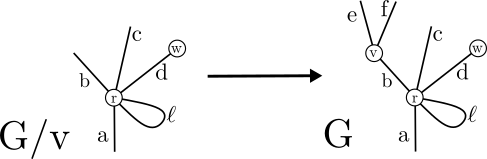}
\end{center}
As a graphical map, the graph $G/v$ is inserted in place of a vertex in an appropriate barbell graph.  This is the opposite way around to inner coface maps when pictured:
\begin{center}
	\includegraphics[scale=0.8]{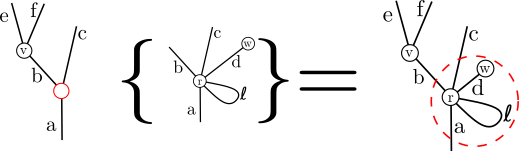}
\end{center}
but it can still be defined as a graphical map, with $\delta^v_0$ the identity on edges and $\delta^v_1$ the identity on vertices.  That is, $\delta^v_1$ sends each vertex in $G/v$ to its associated corolla $C_v$ as an embedding in $G$.

\begin{rmk}
	One important thing to note is that outer vertices are those which have at most one non-leg edge connected to them.
\end{rmk}

\begin{defn}[Cosnip]
\label{def:graphicalCosnip}
	Note that this is a type of outer coface map.
	
	Let $G$ be a graph with a loop $\ell$ (i.e. an edge from a vertex to itself), and define $G/\ell$ to be the graph where $\ell$ has been ``snipped'' to form two edges.  
	Then $\delta^\ell: G/\ell \to G$ is the associated cosnip map
\end{defn}
\begin{center}
	\includegraphics[scale=0.8]{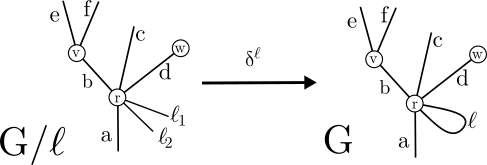}
\end{center}
And as a graphical map, $\delta^\ell_0$ would be surjective on edges, with $\delta^\ell_0(\ell_1) = \delta^\ell_0(\ell_2) = \ell$ and $\delta^\ell_1$ the identity on vertices.  
\begin{center}
	\includegraphics[scale=0.8]{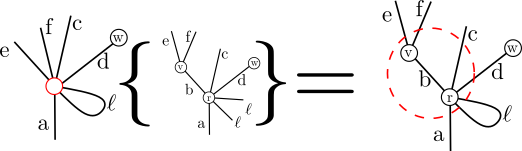}
\end{center}

\begin{defn}[Codegeneracy map]
\label{def:graphicalCodegeneracy}
	Let $b$ be any edge of $G$.  Then define $G_b$ to be the graph where the edge $b$ has been subdivided with exactly one vertex.  Then $\sigma^b: G_b \to G$ is the associated degeneracy map.
\end{defn}
\begin{center}
	\includegraphics[scale=0.8]{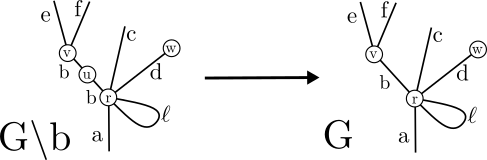}
\end{center}

And as a graphical map $\sigma^b_0$ maps edges to edges according to the labelling, noting that $b$ in $G$ is mapped to by two edges, while $\sigma^b_1$ is given by 
\begin{align*}
	r & \mapsto C_r \\
	u & \mapsto \eta_b \\
	v & \mapsto C_v \\
	w & \mapsto C_w \\
\end{align*}
\begin{center}
	\includegraphics[scale=0.8]{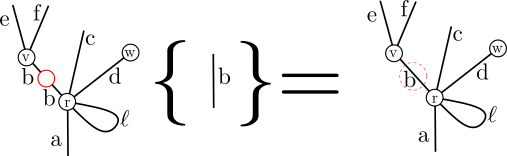}
\end{center}
As is evident above, any of these elementary maps can be written in terms of graphical maps.  The other direction is proven in Lemma \ref{lem:graphtofacedegen}, after the relations are given.

\subsection{Relations between Elementary Graphical Maps}

To make working with this category easier, it is useful to have some relations between these face and degeneracy maps, especially for comparison with other graphical categories like $\cat{\Delta}$.  

\paragraph{Inner coface maps}

Given two (distinct) inner face maps $\delta^a$ and $\delta^b$, $(G/a)/b = (G/b)/a$ and the following diagram commutes:
\begin{center}
\begin{tikzcd}
	(G/b)/a \arrow{r}{\delta^b} \arrow{d}{\delta^a} & G/a \arrow{d}{\delta^a} \\
	G/b \arrow{r}{\delta^b} & G
\end{tikzcd}
\end{center}
In the language of Definition \ref{def:graphicalmap}, both $\delta^a_0 \circ \delta^b_0$ and $\delta^b_0 \circ \delta^a_0$ will be the same inclusion on edges, with $a$ and $b$ not in the image.  Let $a_1$ and $a_2$ denote the vertices on each end of the edge $a$, and likewise let $b_1$ and $b_2$ be the vertices either end of $b$.  Then both $\delta^a_0 \circ \delta^b_0$ and $\delta^b_0 \circ \delta^a_0$ will be the map $V(G/a/b) \to Emb(G)$ that sends the composite vertices $a_1a_2$ and $b_1b_2$ to their associated barbell graphs.  

It is easy to see that this works even in the higher genus case, with loops around vertices,  or with cycles of length greater than or equal to $2$.
\begin{center}
	\label{eg:innerfacerel}
	\includegraphics[scale=0.8]{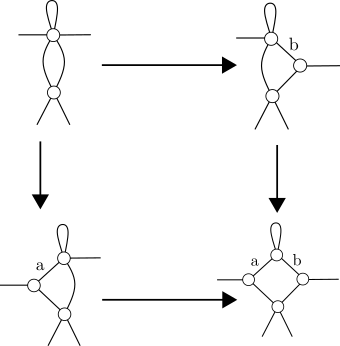}
\end{center}
\begin{center}
	\label{eg:innerfacerel3cyc}
	\includegraphics[scale=0.8]{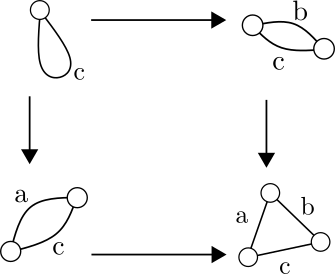}
\end{center}
In the case of a cycle of length $2$, that is, two edges between the same two vertices, once one is contracted it becomes impossible to contract the other.  Thus, if $a$ and $b$ are both edges between vertices $u$ and $v$,  $\delta^a = \delta^b$.

\paragraph{Outer coface maps}
Let $G$ be a tree with at least $3$ vertices, and consider two distinct outer face maps $\delta^v$ and $\delta^w$.  Then $(G/v)/w = (G/w)/v$ and the diagram commutes:
\begin{center}
\begin{tikzcd}
	(T/v)/w \arrow{r}{\delta^v} \arrow{d}{\delta^w} & T/w \arrow{d}{\delta^w} \\
	T/v \arrow{r}{\delta^v} & T
\end{tikzcd}

	\label{eg:RelOuterEG}
	\includegraphics[scale=0.8]{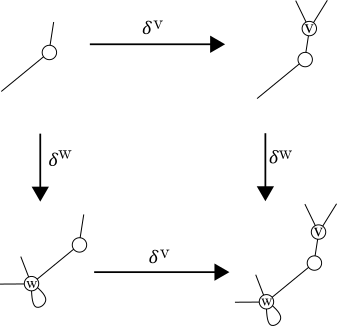}
\end{center}
As with inner coface maps, both $\delta^v \circ \delta^w$ and $\delta^w \circ \delta_v$ will result in the same graphical map, where $(\delta^v \circ \delta^w)_0$ includes on edges (leaving out those outer edges attached to $v$ and $w$) and $(\delta^v \circ \delta^w)_1$ sends each vertex to its associated corolla embedding (essentially, the identity on vertices).

If $G$ is a graph with only two vertices, $v$ and $w$, this still works provided that there is only one edge between them.
\begin{center}
\begin{tikzcd}
	\eta \arrow{r}{\delta^v} \arrow{d}{\delta^w} & G/w \arrow{d}{\delta^w} \\
	G/v \arrow{r}{\delta^v} & G
\end{tikzcd}
\end{center}
\begin{center}
	\label{eg:RelOuterEgTwo}
	\includegraphics[scale=0.8]{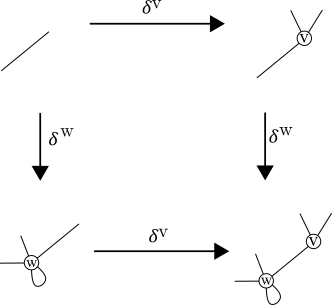}
\end{center}

\paragraph{Inner and Outer}
Then there is the case where an inner coface map is combined with an outer coface map.  Let $\delta^v$ be an outer coface map and $\delta^e$ be an inner coface map.  
If $v$ and $e$ are not adjacent, then $(G/v)/e = (G/e)/v$ and the diagram commutes:
\begin{center}
\begin{tikzcd}
	(G/v)/e \arrow{r}{\delta^e} \arrow{d}{\delta^v} & G/v \arrow{d}{\delta^v} \\
	G/e \arrow{r}{\delta^e} & G
\end{tikzcd}
\end{center}
\begin{center}
	\includegraphics[scale=0.8]{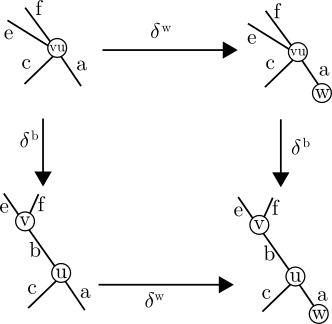}
\end{center}

Suppose $v$ and $e$ are adjacent.  Then denote the vertex on the other side of $e$ by $w$.  The tree $G/e$ combines $v$ and $w$ into a single vertex, denote this $u$.  Then $(G/e)/u$ exists if and only if $(G/v)/w$ exists, $(G/e)/u = (G/v)/w$, and the following diagram commutes:
\begin{center}
\begin{tikzcd}
	(G/e)/u \arrow{r}{\delta^u} \arrow{d}{\delta^w} & G/e \arrow{d}{\delta^e} \\
	G/v \arrow{r}{\delta^v} & G
\end{tikzcd}
\end{center}
\begin{center}
	\includegraphics[scale=0.8]{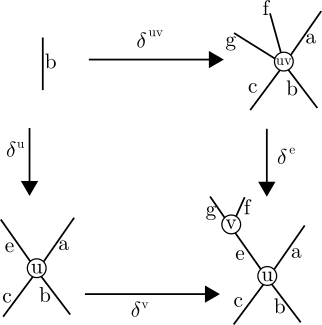}
\end{center}

\paragraph{Degeneracies}

Now consider two codegeneracy maps of $G$ $\sigma^e$ and $\sigma^a$.  Then $(G_e)_a = (G_a)_e $ and the following diagram commutes
\begin{center}
\begin{tikzcd}
	(G_e)_a \arrow{r}{\sigma^e} \arrow{d}{\sigma^a} & G_a \arrow{d}{\sigma^a}\\
	G_e \arrow{r}{\sigma^e} & G
\end{tikzcd}
\end{center}
This describes the map that sends $e_1$ and $e_2$ to $e$, $a_1$ and $a_2$ to $a$, and all other edges to themselves.  All vertices are sent to their associated corollae, except for those extra ones 
which are sent to $\eta$.
\begin{center}
	\includegraphics[scale=0.8]{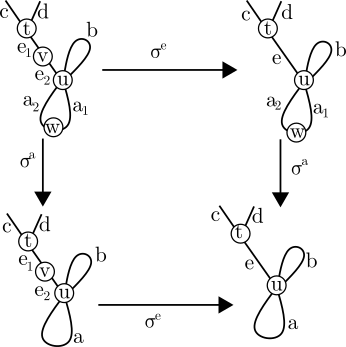}
\end{center}

\paragraph{Face and Degeneracy}
Let $\sigma^a : T_a \to T$ be a codegeneracy and $\delta^e: T' \to T$ a coface map that does not eliminate $a$.   
Then if $\delta^e: T'_a \to T_a$ is the induced coface map, the following diagram commutes:
\begin{center}
\begin{tikzcd}
	T'_e \arrow{r}{\delta^e} \arrow{d}{\sigma^a} & T_e \arrow{d}{\sigma^a} \\
	T' \arrow{r}{\delta^e} & T
\end{tikzcd}
\end{center}
\begin{center}
	\includegraphics[scale=0.8]{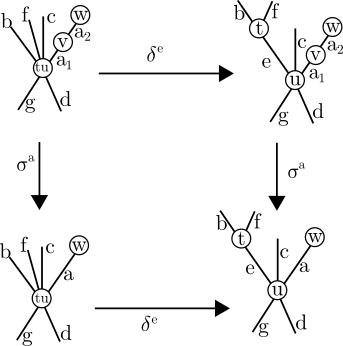}
\end{center}

Consider a codegeneracy map $\delta^e$.  Let the resulting vertex be denoted $v_e$ and the resulting edges be $x$ and $y$.  If $\delta: T \to T^e$ is either
\begin{itemize}
	\item an inner coface map $\delta /x$ or $\delta /y$
	\item an outer coface map $\delta /v_e$
\end{itemize} 
then
\begin{center}
\begin{tikzcd}
	T \arrow{r}{\delta} & T_e \arrow{r}{\sigma^e} & T
\end{tikzcd}
\end{center}
is the identity.

\paragraph{Snip}
Although the cosnip map is a type of outer face map, it acts sufficiently differently to deserve its own treatment.  Firstly, two cosnips commute, as does a snip with an inner coface map.  Both of these compositions can be easily written down as a graphical map, as was done in the previous paragraphs.
Given an inner edge $e$ and a loop $\ell$, 
\begin{center}
\begin{tikzcd}\\
	(G/\ell)/e \arrow{r}{\delta^e} \arrow{d}{\delta^\ell} & G/\ell \arrow{d}{\delta^\ell} \\
	G/e \arrow{r}{\delta^e} & G
\end{tikzcd}
\end{center}
The only interesting part is when the loop is as a result of an inner coface map.  In that case $\ell_1$ and $\ell_2$ are both mapped to $\ell$, and the vertex resulting from the inner face map $e$ is mapped to the barbell containing $e$ as an inner edge, and all other vertices and edges are mapped to themselves. 
\begin{center}
\begin{tikzcd}\\
	(G/\ell)/e \arrow{r}{\delta^e} \arrow{d}{\delta^\ell} & G/\ell \arrow{d}{f} \\
	G/e \arrow{r}{\delta^e} & G
\end{tikzcd}
\end{center}
\begin{center}
	\includegraphics[scale=0.8]{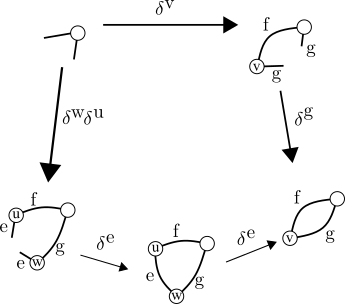}
\end{center}

Given an outer vertex and a disjoint loop, the outer coface map commutes with the snip.  However, when the loop is attached to the outer vertex it becomes a little more complicated.
Let $v$ be an outer vertex containing a loop $\ell$.  Then 
\begin{center}
\begin{tikzcd}
	G/v \arrow{r}{\delta^v} \arrow{dr}{\delta^v} & G/\ell \arrow{d}{\delta^\ell} \\
	& G \\
\end{tikzcd}
\end{center}
This graphical map is the inclusion on both the edges and the vertices, so even though there is a difference between $G/\ell$ and $G$, this difference lies on the vertex $v$ and its associated outer edges, none of which are in the image of this graphical map.

\begin{center}
	\label{eg:outersnipreleg}
	\includegraphics[scale=0.8]{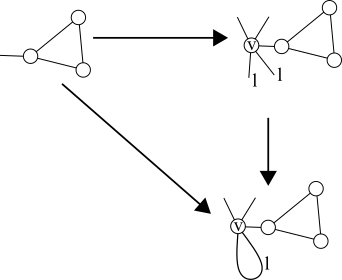}
\end{center}
There are also a couple of relations between the snip and an outer face map.  Let $v$ be an outer vertex connected to two inner edges $f$ and $g$.   If $v$ is split into two vertices by an inner coface map, let these be labelled $u$ and $w$, and the relevant inner edge be labelled $e$.  Then 
	$$\delta^g \delta^v = \delta^e \delta^{e'} \delta^u \delta^w.$$
	Where, to disambiguate, the $\delta^{e'} $ refers to cosnipping the edge $e$ and $\delta^e$ refers to contracting the edge $e$.
\begin{center}
	\label{eg:outersniprelalteg}
	\includegraphics[scale=0.8]{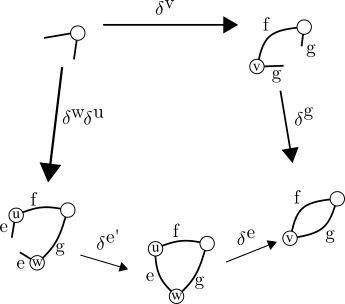}
\end{center}
This is because $\delta^g$ and $\delta^e \delta^{e'}$ only affect those edges and vertices which are not in the image of $\delta^v$ and $\delta^w\delta^u$.

Then, one may wonder whether there are any tricks related to the following image, but as is clear below this only occurs when certain edges are the same colour, and these maps commute anyway.  Again, the composite map is the inclusion, and the differences are happening on a vertex that is not in the image.

\begin{center}
	\label{eg:outersniprelinnereg}
	\includegraphics[scale=0.8]{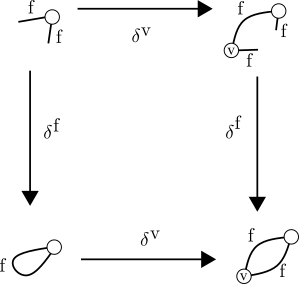}
\end{center}

For ease of access, these relations will once again be arranged in a table.  The heading on the right indicates which map is done first, and then the left column says which map is done second.  Those relations between vertices and edges which are disjoint will be omitted, as they always commute, so in this table assume the face and degeneracy maps are always coincident.  Additionally, the elementary maps labelled in the columns are applied first, followed by the rows.
	\begin{center}
	\begin{tabular}{|c|c|}
	\hline 
	 & inner \\  
	\hline 
	inner & $\delta^{\overline{pq}}\delta^{\overline{qr}}=\delta^{\overline{qr}}\delta^{\overline{pq}}$ \\ 
	\hline 
	outer & $\delta^{pq}\delta^{\overline{pq}} = \delta^q \delta^p$ \\ 
	\hline 
	degeneracy & $\sigma^{\overline{p(qr)}} \delta^{\overline{qr}} = id = \delta^{\overline{qr}} \sigma^{\overline{pq}}$ \\
	\hline
	snip & $\delta^{\overline{pq}} \delta^{\overline{(pq)(pq)}} = f \delta^{\overline{pq}}$ \\
	\hline
	& outer \\
	\hline 
	inner & $\delta^{\overline{pq}}\delta^p$ undefined \\ 
	\hline 
	outer & $\delta^q \delta^p = \delta^{pq} \delta^{\overline{pq}}$ \\ 
	\hline 
	degeneracy & $\sigma^p \delta^q = \delta^q \sigma^{\overline{pq}}$ \\
	\hline
	snip & $\delta^{\overline{pp}} \delta^p = \delta^p$  \\
	\hline
	& degeneracy \\
	\hline
	inner & $\delta^{\overline{p (pq)}} \sigma^{\overline{pq}} = id = \sigma^{\overline{r(pq)}} \delta^{\overline{pq}}$ \\
	\hline
	outer & $\sigma^p \delta^q = \delta^q \sigma^{\overline{pq}}$ \\
	\hline
	degeneracy & $\sigma^{\overline{p(pq)}} \sigma^{\overline{pq}} = \sigma^{\overline{(pq)q}} \sigma^{\overline{pq}}$ \\
	\hline
	snip & undefined or commutes \\
	\hline
	\end{tabular} 
	\end{center}
	If an edge has been removed by face map, then to snip it is undefined.  If there is a vertex of degree 2, and one incident edge is snipped, then the other map cannot be, for this would make the graph disconnected.  Otherwise, the snip commutes with itself.  In addition, there is another relation between a snip coincident to an outer face map.  If an outer vertex $v$ is connected to two inner edges $f$ and $g$, and $v$ is split into two vertices by an inner coface map, let these be labelled $u$ and $w$, and the relevant inner edge be labelled $e$.  Then 
	$$\delta^f \delta^v = \delta^e \delta^{e'} \delta^u \delta^w.$$

\begin{lem}
	\label{lem:graphrelationsmovingpast}
	Given any sequence of coface and codegeneracy maps, one can find a standard form consisting of a sequence of codegeneracy maps followed by a sequence of coface maps.
\end{lem}
\begin{proof}
	Consider the above table of relations.  If a codegeneracy is coincident to an inner coface map they will annihilate.  All other maps commute with codegeneracy maps.  Therefore they can be moved past each other to form a sequence of codegeneracy maps followed by coface maps.
\end{proof}

One important thing to note is that morphisms, when written in terms of coface and codegeneracy maps, are slightly different in form to morphisms written in terms of graphical maps.  The former usually take the form of a sequence of maps $f^a : G' \to G$, where $a$ represents the vertex or edge to be removed or added.  The latter take the form $G' = G{H_v}$ or $G' = H{G_v}$, where the $H$ contains the information about which vertex or edge is to be modified.  Theorem \ref{lem:graphtofacedegen} gives a proof of their equivalence.  An alternative proof of this theorem can also be found in \cite[Theorem 2.7]{robertson2019modular}.

\begin{thm}
	\label{lem:graphtofacedegen}
	Any morphism of graphs $G \to G''$ can be decomposed into a composition of codegeneracy maps followed by a composition of coface maps. 
\end{thm}
\begin{proof}
	By Lemma \ref{lem:graphicalmapgraphsub} 
	each graphical map can be written as a series of graph substitutions followed by an inclusion.  Denote this $G \to G' \to G''$.  By Lemma \ref{lem:graphrelationsmovingpast} a series of coface and codegeneracy maps can be converted to a standard form, so it suffices to show that each graph substitution and inclusion can be decomposed into a series of coface and codegenercy maps.
	
	We will proceed by induction on the number of vertices in $G''$.  Firstly, consider a graphical map $G \to \eta$.  The only inclusion is $\eta \inj \eta$ itself, and the only graph substitutions available are codegeneracies.  This can be written as a series of graph substitutions followed by an inclusion. 

	  Consider a single graph substitution into a vertex v, $G_v' = G\{H_v\}$, where $H_v$ has fewer vertices than $G''$.  By the inductive hypothesis, we can write $H_v$ as a series of inner coface and codegeneracy maps applied to the corolla $C_v$.  Then, because these inner coface and codegeneracy maps only involve that part of $G'$ which is contained in $H_v$, we can move them outside and produce the following.
	\begin{align*}
		G_v' 
		& = G\{ H_v \} \\
		& = G\{ \delta^1 \ldots \delta^m \sigma^1 \ldots \sigma^n C_v\} \\
		& = \delta^1 \ldots \delta^m \sigma^1 \ldots \sigma^n G\{  C_v\} \\
		& = \delta^1 \ldots \delta^m \sigma^1 \ldots \sigma^n G \\
	\end{align*}
	This process can be repeated for each subsequent vertex and associated graph substitution, with the composition of these maps being a map $G \to G'$

	Then, there is the inclusion $G' \to G''$.  This can be written as a graph substitution $G'' = H\{G'\}$, where $H$ is the graph obtained by replacing $G'$ in $G''$ with its associated corolla.  That is, the corolla which is identical on its boundary.  Again, the graph $H$ can be built up from a series of coface and codegeneracy maps, starting from $C_v$ and without interfering with that vertex.  Then 
	\begin{align*}
		G'
		& = H\{ G_v \} \\
		& = \delta^1 \ldots \delta^m \sigma^1 \ldots \sigma^n C_v \{ G \} \\
		& = \delta^1 \ldots \delta^m \sigma^1 \ldots \sigma^n G \\
	\end{align*}
\end{proof}	
\end{appendices}

\bibliographystyle{plainnat}
\bibliography{phd}

\end{document}